\documentclass[journal,twoside,web]{ieeecolor}
\usepackage{generic}
\usepackage{cite}
\usepackage{amsmath,amssymb,amsfonts}
\usepackage{algorithmic}
\usepackage{graphicx}
\usepackage{algorithm,algorithmic}
\usepackage{hyperref}
\hypersetup{hidelinks=true}
\usepackage{textcomp}

\usepackage{mathdots}
\usepackage{mathbbol}

\def\BibTeX{{\rm B\kern-.05em{\sc i\kern-.025em b}\kern-.08em
    T\kern-.1667em\lower.7ex\hbox{E}\kern-.125emX}}
\begin{document}
\title{Spectral Decompositions of Controllability Gramian and Its Inverse based on System Eigenvalues in Companion Form}
\author{Iskakov Alexey and Yadykin Igor
\thanks{This work was supported by the Russian Science Foundation, project no. 25-29-20158.}
\thanks{Iskakov A.B. is with the V. A. Trapeznikov Institute of Control Sciences of Russian Academy of Sciences, 
65 Profsoyuznaya street, Moscow 117997, Russia (e-mail: iskalexey@gmail.com).}
\thanks{Yadykin I.B. is with the V. A. Trapeznikov Institute of Control Sciences of Russian Academy of Sciences, 
65 Profsoyuznaya street, Moscow 117997, Russia (e-mail: jad@ipu.ru).}}

\maketitle

\begin{abstract}
Controllability and observability Gramians, along with their inverses, are widely used to solve various problems in control theory. This paper proposes spectral decompositions of the controllability Gramian and its inverse based on system eigenvalues for a continuous LTI dynamical system in the controllability canonical (companion) form. The Gramian and its inverse are represented as sums of Hermitian matrices, each corresponding to individual system eigenvalues or their pairwise combinations. These decompositions are obtained for the solutions of both algebraic and differential Lyapunov and Riccati equations with arbitrary initial conditions, allowing for the estimation of system spectral properties over an arbitrary time interval and their prediction at future moments. The derived decompositions are also generalized to the case of multiple eigenvalues in the dynamics matrix spectrum, enabling a closed-form estimation of the effects of resonant interactions with the system's eigenmodes. The spectral components are interpreted as measurable quantities in the minimum energy control problem. Therefore, they are unambiguously defined and can quantitatively characterize the influence of individual eigenmodes and associated system devices on controllability, observability, and the asymptotic dynamics of perturbation energy. The additional information obtained from these decompositions can improve the accuracy of algorithms in solving various practical problems, such as stability analysis, minimum energy control, structural design, tuning regulators, optimal placement of actuators and sensors, network analysis, and model order reduction.
\end{abstract}
\section{Introduction}
{\it Controllability and observability Gramians} play an important role in the qualitative and quantitative evaluation of controllability, observability, and stability for different classes of dynamical systems. For LTI systems, the observability Gramian determines the integral measure of the output signal with zero control. Conversely, the inverse of the controllability Gramian determines the minimum energy required to bring the system to a reachable state. The concept of Gramians has been generalized and interpreted for deterministic bilinear, stochastic linear, and switched linear systems \cite{Benner-2011}, \cite{Duff-2019}. Gramians and their inverses are widely used to solve various practical problems in control theory, such as structural design \cite{Gupta-2020}, minimum energy control \cite{Casadei-2020}, model reduction \cite{Benner-2021}, optimal selection of actuators and sensors \cite{Dilip-2019}, network analysis and design \cite{Zhao-2019}, analysis of oscillations and their interactions \cite{Isk-2023}, tuning of controllers \cite{Ghosh-2020}, optimal placement of phasor measurement units in power grids \cite{Qi-2015}, and control of traffic networks \cite{Bianchin-2020}.

{\it The spectral properties of Gramians} are effectively used in model order reduction methods based on the singular value decomposition of Gramians (see the review in \cite{Baur-2014}). These methods find applications in Gramian-based filtering schemes \cite{SVD-Filter-2022}, image processing \cite{SVD-Image-2023}, and approximation of large sparse power systems \cite{SVD-Pow-Sys-2008}. An alternative approach, where Gramians are decomposed not over their own spectrum but over the spectrum of the dynamics matrix of the system was considered in \cite{Antoulas-2005} and in a more general form in \cite{Yad-2010}, \cite{DAN-2017}, \cite{Behr-2019}, and \cite{Yad-2022}. The Gramian is represented as a sum of matrices corresponding to individual system eigenvalues or their pairwise combinations, which allows for a quantitative characterization of the influence of individual eigenmodes and associated system devices on controllability, observability, and the asymptotic dynamics of perturbation energy. We believe that these Spectral Decompositions based on System Eigenvalues (SDSE) can combine the advantages of modal analysis of the system with Lyapunov stability analysis and
improve various existing methods for stability estimation and optimal control, which are based on the application of Gramians \cite{Automatica-2021}.
In particular, in \cite{Isk-2023} SDSE of Gramians has been applied to reveal the structure and resonance interaction of inter-area oscillations in power systems. 

This paper proposes analytical expressions for SDSE of the controllability Gramian and its inverse for a continuous LTI dynamical system. Due to the duality of controllability and observability, the spectral decompositions for the observability Gramian are obtained similarly and are not considered separately. {\it Controllability and observability canonical forms},  characterized by the minimum number of non-zero parameters describing system dynamics, have proven to be a very convenient framework for solving control design and observer design problems. In this paper, we also consider Gramians in the canonical form of controllability or observability to obtain their spectral decompositions in the simplest and most concise form, which does not depend on input and output matrices.
Although SDSE of Gramians have been studied previously in \cite{Antoulas-2005}, \cite{Hauksdottir-2008}, \cite{Yad-2010}, \cite{DAN-2017}, \cite{Behr-2019}, and \cite{Yad-2022}, this paper combines and extends these results. Specifically, we derive analytic expressions for the spectral decomposition of the solution to the differential Lyapunov equation with an arbitrary initial condition (in Theorem 2 and Remark 3), allowing for the estimation of system spectral properties over an arbitrary time interval and their prediction at future moments. The derived decompositions are also generalized to the case of multiple eigenvalues in the dynamics matrix spectrum (in Theorems 6 and 7), enabling a closed-form estimation of the effects of resonant interactions with the system's eigenmodes. 
To the best of our knowledge, spectral decomposition for the inverse of the Gramian is obtained for both the algebraic and differential Riccati equations for the first time (Theorems 4, 5, and 7). An additional result is a closed-form solution of a certain class of Riccati equations (in Corollary 1).

The paper is organized as follows: Section II presents the problem statement and provides preliminary results on the spectral decomposition of Gramians from the literature. In Section III, we derive SDSE for the finite controllability Gramian with an arbitrary initial condition and provide an interpretation for its spectral components. Section IV shows how to obtain SDSE of the Gramian for an arbitrary controllable multiple-input (MI) system from the previously obtained SDSE in the canonical controllability form.
In Section V, SDSE of the inverse of Gramian in the canonical controllability form is derived for both algebraic and differential Riccati equations. 
Section VI extends the obtained results to the case of multiple eigenvalues in the dynamics matrix spectrum. Section VII concludes the paper.\\

\section{Problem Formulation and Preliminaries}
Consider an {\it LTI dynamical system} described by the equations
\begin{equation}  \label{sys-0}
\dot{x}(t) =A x(t) + B u(t), 
\end{equation}
where $x(t)\in \mathbb{R}^n$ is the state of the system, $u(t)\in \mathbb{R}^m$ is the input vector, and the constant matrices $A\in \mathbb{R}^{(n\times n)} and \ B\in \mathbb{R}^{(n \times m)}$ describe system dynamics and input, respectively. We assume that $u(t)$ satisfies conditions that guarantee the existence of an absolutely continuous solution. Throughout this paper, except in Section VI, we assume that matrix $A$ is diagonalizable and has a complete set of eigenvectors.

The matrix {\it Lyapunov algebraic and differential equations} for system \eqref{sys-0} are
\begin{equation}  \label{Lyap-algeb} 
 A P+ P A^T= - BB^T ,                                      
\end{equation}
\begin{equation}  \label{Lyap-diff} 
- \frac{dP(t)}{dt} + A P(t) + P(t) A^T= - BB^T .                                     
\end{equation}
They have a unique solution if and only if all eigenvalues from the spectrum of matrix $A$ satisfy the condition
\begin{equation} \label{Lyap-cond}
\lambda_i+\lambda_j \ne 0 \ \ \text{for all} \ \  \lambda_i,\lambda_j\in \sigma(A) , \ \ i,j = \overline{1,n} , 
\end{equation}
which we assume to hold. For a strictly stable $A$, the solutions of \eqref{Lyap-algeb} and \eqref{Lyap-diff} can be interpreted 
as {\it the infinite and finite controllability Gramians}, respectively. 
Note that the SDSEs of the solution of \eqref{Lyap-algeb} and \eqref{Lyap-diff} obtained in this paper will also be valid 
for the unstable matrix $A$.

Let matrix $A$ have {\it the characteristic equation} $\mathrm{det}(Is - A) = \sum^n_{i=0} a_i s^i$, where $a_n=1$, and let $\mathcal{C} = [ B,AB,A^2 B,...,A^{(n-1)} B]$ be {\it the controllability matrix}. Consider a single-input (SI) controllable system where $m=1$ and $B=b\in R^n$. 
Then, according to \cite{Hauksdottir-2009}, there exists a non-singular similarity matrix 
\begin{equation} \label{transfer-matrix}
T = \mathcal{C} \mathcal{H}_u, \ \ \mathcal{H}_u =      
\begin{bmatrix}
a_1 & a_2 & \cdots & a_{n-1}  & 1  \\ 
a_2 & \iddots & \iddots  & 1            & 0  \\
\vdots &  \iddots   & \iddots & \iddots & \vdots \\
a_{n-1} & 1 & 0    &  \iddots      & 0  \\
1        & 0   &  \cdots      &  0     & 0  
\end{bmatrix} ,
\end{equation}
that transforms \eqref{sys-0} into {\it the controllability canonical form}
\begin{gather} 
\dot{x}_c (t) = A_C \, x_c(t) + b_C \, u(t) , \ \ \text{where} \ \nonumber \\              
x = T x_c, \  A \, T = T A_C, \  B = T b_C, \label{sys-c} \\   
A_C = 
\begin{bmatrix}
0_{(n-1)\times1} &  & I_{(n-1)\times(n-1)} &   \\ 
-a_0 & -a_1 & \dots & -a_{n-1}              
\end{bmatrix} ,  \nonumber \\
b_C= [ 0, 0, \cdots, 0, 1]^T. \nonumber
\end{gather}
Lyapunov equations \eqref{Lyap-algeb} and \eqref{Lyap-diff} are transformed to 
\begin{gather}   
 A_C P_C + P_C A_C^T = - \, b_C b_C^T ,  \label{Lyap-algeb-c} \\ 
-\frac{dP_C(t)}{dt} + A_C P_C(t) + P_C(t) A_C^T = - \, b_C b_C^T ,   \label{Lyap-diff-c} \\ 
P = T P_C \, T^T .    \label{P-to-orig-sys}  
\end{gather}
The infinite controllability Gramian $P_C$ for system \eqref{Lyap-algeb-c} in the canonical form has 
a particularly simple zero-plaid Hankel alternating structure 
\cite{Sreeram-1991,Xiao-1992}:
\begin{equation} \label{zero-plaid}     
P_C =     
\begin{bmatrix}
p_1 & 0     & -p_2 & 0  & p_3 & \dots  \\ 
0     & p_2 & 0    & -p_3 & 0 & \dots    \\
-p_2 &  0  & p_3 &   0  & -p_4 & \dots \\
 0    & -p_3 & 0    &p_4 & 0     & \dots  \\
p_3  & 0    & -p_4 &  0 &\ddots & \dots \\
 \dots &  \dots &  \dots &  \dots &  \dots & p_n 
\end{bmatrix} \ .
\end{equation}
In \cite{Hauksdottir-2009}, an algorithm is proposed to calculate the $p_i$ coefficients of $P_C$ in \eqref{zero-plaid} from the $a_i$-coefficients of the characteristic polynomial using O($n^2$) operations. 

To  find the controllability Gramian $P$ from $P_C$ in the general case of an MI system, where $m > 1$, we can efficiently treat each of the $m$ columns of $B$ separately using the expression proposed in \cite{Hauksdottir-2009}:
\begin{equation}  \label{MI-P-from-Pc}  
P = \mathcal{C} \, (\mathcal{H}_u P_C \mathcal{H}_u  \otimes I_m) \, \mathcal{C}^* ,    
\end{equation}
where $I_m$ is a unit matrix $m\times m$, $(\cdot)^*$ is the complex conjugate transpose, and $\otimes$ is the matrix Kronecker product. \\

In the subsequent presentation, we will use the following result as a starting point. The finite controllability Gramian $P(t)$ can be obtained as a solution of the differential Lyapunov equation \eqref{Lyap-diff}
 and can be represented as a sum of Hermitian matrices corresponding to either individual eigenvalues or their pairwise combinations 
 \cite{Yad-2010, DAN-2017}:
 \begin{equation}  \label{SD-gen-diff}
  P(t) = \sum_{i} \tilde{P}_i(t) = \sum_{i,j} P_{ij}(t) , \ \ \tilde{P}_i(t) = \sum_{j} P_{ij}(t) .    
\end{equation}
All summation indices here and elsewhere run from $1$ to $n$, unless otherwise specified. If the spectrum of the dynamics matrix is simple, $\sigma(A)=\{\lambda_1,\lambda_2,\cdots,\lambda_n\}$, then SDSE of the finite controllability Gramian with zero initial condition is computed as
\begin{gather} 
\tilde{P}_i(t) = -\left\{R_i B B^T (\lambda_i I + A^*)^{-1}(I- e^{(I\lambda_i+A^*)t})\right\}_H \, , \nonumber \\
P_{ij}(t) =\left\{\frac{e^{(\lambda_i+\lambda^*_j)t} -1}{\lambda_i+\lambda^*_j} \cdot R_i B B^T R^*_j\right\}_H ,   \label{subgram-general}    
\end{gather}
where  
$\{\cdot\}_H$ represents the Hermitian part of a matrix, and $R_i, R_j$ are the matrix residues, that is, the coefficients in the expansion of the resolvent of matrix $A$ corresponding to its eigenvalues $\lambda_i$ and $\lambda_j$:
 \begin{equation}  \label{simple-residues}
 (Is-A)^{-1} = \frac{R_1}{s-\lambda_1} + \frac{R_2}{s-\lambda_2} + \cdots + \frac{R_n}{s-\lambda_n} . \nonumber
\end{equation}
The spectral components $\tilde{P}_i(t)$ and $P_{ij}(t)$ in \eqref{subgram-general} were named {\it sub-Gramians} and {\it pair sub-Gramians}, respectively. \\

\section{SDSE of Gramians in Companion Form}

In this section, Theorem 1 first considers the SDSE of the infinite Gramian in the controllability canonical form 
obtained in \cite{Yad-2022} using the Faddeev-Leverrier algorithm.
We will show that the spectral components in this decomposition are uniquely defined as measurable quantities in the minimum energy control problem. To analyze the spectral properties of the system on an arbitrary time interval and predict their dynamics at future times, we extend the result of Theorem 1 to the case of a finite Gramian with an arbitrary initial condition in Theorem 2 and Remark 3. \\

\noindent \textbf{Theorem 1} \cite{Yad-2022}. \textit{
Let system \eqref{sys-c} have a simple spectrum. Then, the SDSE of \textbf{the infinite} controllability Gramian $P_C$ in \eqref{Lyap-algeb-c} is}
\begin{gather}  
  P_C = \sum_{i} \tilde{P}^C_i = \sum_{i,j} P^C_{ij} , \ \ \tilde{P}^C_i = \sum_{j} P^C_{ij} ,   \label{SD-alg-c} \\
\tilde{P}^C_i = \left\{ \frac{x_i x^T_i \mathcal{J}}{N'(\lambda_i) N(-\lambda_i)}\right\}_H \, , \label{subG-alg-c} \\
P^C_{ij} = \left\{ \frac{-1}{\lambda_i+\lambda^*_j} \frac{x_i x^*_j}{N'(\lambda_i) N'(\lambda^*_j)} \right\}_H \, , \label{subG2-alg-c}\\
 x_i = [1,\lambda_i, \lambda^2_i, \cdots, \lambda^{n-1}_i ]^T \ , \nonumber
\end{gather}
\textit{where $N(s)=\text{det}(Is-A)$ and $N'(s)$ are the characteristic polynomial and its derivative, respectively, and $\mathcal{J} = \mathrm{diag}\,(-1, (-1)^2, \dots,(-1)^n)$.}\\

\noindent \textbf{Remark 1.} The element-wise representation of matrices $\tilde{P}^C_i$ and $P^C_{ij}$ in (\ref{subG-alg-c}-\ref{subG2-alg-c}) takes the following form:
\begin{gather} 
\tilde{P}^C_i = \sum_{\mu, \nu} 
\left\{\frac{(-1)^{\nu-1} \lambda_i^{\mu+\nu-2}}{N'(\lambda_i) N(-\lambda_i)} \ \mathbb{1}_{\mu \nu}\right\}_H \, , \nonumber \\
P^C_{ij} =  \sum_{\mu, \nu}
\left\{\frac{ -1}{\lambda_i+\lambda^*_j} \cdot \frac{\lambda_i^{\mu-1} (\lambda^*_j)^{\nu-1}}{N'(\lambda_i) N'(\lambda^*_j)} \mathbb{1}_{\mu \nu} \right\}_H , \label{subG-alg-c-elem}
\end{gather}
where $\mathbb{1}_{\mu \nu} = e_{\mu} e_{\nu}^T$ is the vector product of the $\mu$ column and $\nu$ row of the unit matrix.\\

Note that the stability of system \eqref{sys-c} is not required in Theorem 1 for \eqref{SD-alg-c}, \eqref{subG-alg-c} and \eqref{subG2-alg-c} to hold, an assumption made in \cite{Yad-2022}. This will be evident from the proof of Theorem 2 below. In this case, however, $P_C$ should be interpreted simply as a solution of \eqref{Lyap-algeb-c}, not necessarily as an infinite controllability Gramian.

Let us investigate some properties of the sub-Gramians $\tilde{P}^C_i$ and $P^C_{ij}$ in Theorem 1. They are not merely mathematical entities but also measurable quantities. They can be interpreted using the notion of {\it minimum energy control}. 
If the state $x_0$ is reachable, then the minimum energy $E_{min}$ to bring stable LTI system \eqref{sys-c} from the zero state to $x_0$ and the corresponding optimal control $\hat{u}(t)$ are determined by the inverse of the controllability Gramian \cite{Benner-2011}:
\begin{gather}
E_{min} = \inf_{x(-\infty)=0} \int^0_{-\infty} \hat{u}^T(t) \hat{u}(t) dt = x^T_0 P^{-1}_C x_0 , \ \ \text{where} \nonumber \\
 \hat{u}(t) = e^T_n e^{-A^T_C t} P^{-1}_C x_0 , \ \ -\infty< t< 0 \label{Emin}
\end{gather}
is the optimal control. It is divided into spectral components:
\begin{equation}
 \hat{u}(t) = \sum_i \hat{u}_i (t), \ \ \hat{u}_i (t) = e^T_n R^*_i e^{-\lambda^*_i t} P^{-1}_C x_0 . \label{Umin}
\end{equation}
\noindent \textbf{Property 1 (Interpretation).}  
\textit{Let system \eqref{sys-c} be stable and state $x_0$ be reachable. The controllability sub-Gramian $\tilde{P}^C_i$ in problem \eqref{Emin} characterizes the integral of the product of $\hat{u}(t)$ and its $i$-th modal component. The pair controllability sub-Gramian $P^C_{ij}$ characterizes the integral of the product of the $i$-th and $j$-th modal components of $\hat{u}(t)$, i.e.,}
\begin{gather}
x^T_0 ( P^{-1}_C)^T \tilde{P}^C_i P^{-1}_C x_0 = \frac{1}{2} \int^0_{-\infty} (\hat{u}^*_i \hat{u} + \hat{u}^* \hat{u}_i) dt , \nonumber \\
x^T_0 ( P^{-1}_C)^T P^C_{ij} P^{-1}_C x_0 = \frac{1}{2} \int^0_{-\infty} (\hat{u}^*_i \hat{u}_j + \hat{u}^*_j \hat{u}_i) dt . \nonumber
\end{gather}

\noindent \textbf{The proof} is obtained by directly substituting \eqref{subG-alg-c}, \eqref{subG2-alg-c}, \eqref{Emin} and \eqref{Umin} into these equations and verifying them.\\ 

This interpretation of controllability sub-Gramians carries over to an arbitrary controllable system \eqref{sys-0} after transitioning to the corresponding coordinate system. Therefore, the sub-Gramians are uniquely defined as measurable quantities in problem \eqref{Emin}. A similar interpretation for observability sub-Gramians is given in \cite{Isk-2023} (Theorem 1) for disturbances in power systems. A similar energy-based interpretation of the spectral components of the squared $H_2$ norm of the transfer function is proposed in \cite{Hauksdottir-2008}.\\

\noindent \textbf{Property 2 (Structure).} \textit{
The sub-Gramians $\tilde{P}^C_i$ preserve the zero-plaid Hankel alternating structure \eqref{zero-plaid}, i.e.
their elements satisfy the formulas}
\begin{equation}  \label{exam2}
(\tilde{P}^C_i)_{\mu \nu} = \left\{
                \begin{array}{ll}
                (-1)^{\nu-1} (\tilde{P}^C_i)_{m m}, \ \mu+\nu = 2m \\
                0 , \ \ \mu+\nu = 2m+1
                \end{array}                
\right., \ m = \overline{1,n} \, .  \nonumber 
\end{equation}
 
\noindent \textbf{Proof.} Consider in \eqref{subG-alg-c-elem} the expressions 
$(\hat{P}^C_i)_{\mu \nu} = \frac{(-1)^{\nu-1} \lambda_i^{\mu+\nu-2}}{N'(\lambda_i) N(-\lambda_i)}$ 
 under $\{\cdot\}_H$ at positions with an odd index $\mu+\nu =2m+1$. For these, the following holds:
\begin{gather} 
(\hat{P}^C_i)_{\mu \nu} / (-1)^{\nu-1}  = (\hat{P}^C_i)_{\nu \mu} / (-1)^{\mu-1} , \ \  \text{or} \nonumber \\
(\hat{P}^C_i)_{\mu \nu} = (-1)^{\nu-\mu} (\hat{P}^C_i)_{\nu \mu} \ \ \ \ \ \ \ \ \ \ \ \ \ \ \ \ \ \ \ \ \ \  \ \nonumber \\
\ \ \ \ \ \ \ \ \ \ =(-1)^{1+2m-2\mu}(\hat{P}^C_i)_{\nu \mu} = - (\hat{P}^C_i)_{\nu \mu} . \nonumber
\end{gather}
Therefore, after the Hermitian conjugation operation, the corresponding elements of matrix $\tilde{P}^C_i$ will become zero. The sign alternation of the elements along the diagonals defined by the condition $\mu+\nu = 2m$ is evident. $\blacksquare$\\

It follows from Property 2 and decomposition \eqref{SD-alg-c} that the Gramian $P_C$ preserves the zero-plaid Hankel alternating structure \eqref{zero-plaid}, which is an independent derivation of the result obtained in \cite{Sreeram-1991, Xiao-1992}. The elements of the pair sub-Gramians $P_{ij}^C$ in \eqref{subG-alg-c} no longer have such a Hankel structure, 
but the following property holds for them.\\

\noindent \textbf{Property 3.} \textit{If system \eqref{sys-c} is stable, then 
$P_{ii}^C \ge 0$, $i = \overline{1,n}$ .}\\

The proof follows from the fact that the real parts of the eigenvalues of a stable system are negative, and $P_{ii}^C$ is proportional to $[1,\lambda_i,\lambda_i^2,\dots,\lambda_i^n]^T \cdot ([1,\lambda_i,\lambda_i^2,...,\lambda_i^n]^T)^*\ge 0$.\\
 
\noindent \textbf{Example 1} illustrates an application of Theorem 1. 
Consider system \eqref{sys-c} of dimension $3\times3$ with eigenvalues $\lambda_1=1,\lambda_2=2,\lambda_3=3$. 
Its characteristic polynomial and dynamics matrix take the form
\begin{equation}
N(s)=s^3-6s^2+11s-6 , \ \ A_C = 
\begin{bmatrix}
0 & 1 & 0 \\
0 & 0 & 1 \\
6 & -11 & 6    
\end{bmatrix} . \nonumber
\end{equation}
According to \eqref{subG-alg-c}, we compute the sub-Gramians:
\begin{gather}
\tilde{P}^C_1 = \frac{-1}{48}
\begin{bmatrix}
1 & 0 & 1 \\
0 & -1 & 0 \\
1 & 0 & 1    
\end{bmatrix} , \ \tilde{P}^C_2 = \frac{1}{60}
\begin{bmatrix}
1 & 0 & 4 \\
0 & -4 & 0 \\
4 & 0 & 16   
\end{bmatrix} , \nonumber \\
\tilde{P}^C_3 = \frac{-1}{240}
\begin{bmatrix}
1 & 0 & 9 \\
0 & -9 & 0 \\
9 & 0 & 81    
\end{bmatrix} . \nonumber 
\end{gather}
Their sum gives the solution of the Lyapunov equation \eqref{Lyap-algeb-c}:
\begin{equation}
P_C = \tilde{P}^C_1 + \tilde{P}^C_2 + \tilde{P}^C_3 = \frac{-1}{120}
\begin{bmatrix}
1 & 0 & -1 \\
0 & 1 & 0 \\
-1 & 0 & 11    
\end{bmatrix} . \nonumber
\end{equation}
According to \eqref{subG-alg-c}, we also obtain the pair sub-Gramians:
\begin{gather}
P^C_{11} = \frac{-1}{8}
\begin{bmatrix}
1 & 1 & 1 \\
1 & 1 & 1 \\
1 & 1 & 1    
\end{bmatrix} , \ 
P^C_{12} = P^C_{21} = \frac{1}{12}
\begin{bmatrix}
2 & 3 & 5 \\
3 & 4 & 6 \\
5 & 6 & 8    
\end{bmatrix} , \ \nonumber \\
P^C_{22} = \frac{-1}{4}
\begin{bmatrix}
1 & 2 & 4 \\
2 & 4 & 8 \\
4 & 8 & 16   
\end{bmatrix} , \ 
P^C_{23} = P^C_{32} = \frac{1}{20}
\begin{bmatrix}
2 & 5 & 13 \\
5 & 12 & 30 \\
13 & 30 & 72    
\end{bmatrix} , \nonumber \\
P^C_{33} = \frac{-1}{24}
\begin{bmatrix}
1 & 3 & 9 \\
3 & 9 & 27 \\
9 & 27 & 81    
\end{bmatrix} , \
P^C_{31} = P^C_{13} = \frac{-1}{16}
\begin{bmatrix}
1 & 2 & 5 \\
5 & 3 & 6 \\
5 & 6 & 9    
\end{bmatrix} . \nonumber 
\end{gather}
One can verify that:
\begin{gather}
\tilde{P}^C_1 = P^C_{11} + P^C_{12} + P^C_{13},  \ \tilde{P}^C_2 = P^C_{21} + P^C_{22} + P^C_{23}, \nonumber \\
\tilde{P}^C_3 = P^C_{31} + P^C_{32} + P^C_{33},  \ \ P_C = \sum^n_{i,j=1} P^C_{ij} \ \square \ . \nonumber 
\end{gather}

The SDSE of {\it the infinite} Gramian obtained in Theorem 1 allows analysis of only the asymptotic spectral properties of the system. To study these properties over an arbitrary time interval and predict their dynamics at future times, we derive the SDSE of {\it the finite} controllability Gramian, which is a solution to \eqref{Lyap-diff-c} for an arbitrary initial condition. In contrast to \cite{Yad-2022}, we will carry out the proof using the eigenvectors of the dynamics matrix  
and will not assume system stability. \\ 

\noindent \textbf{Theorem 2.} \textit{
Let system \eqref{sys-c} have a simple spectrum. Then, SDSE of \textbf{the finite} controllability Gramian $P_C(t)$, 
satisfying \eqref{Lyap-diff-c} with zero initial condition $P_C (0)=0$, is represented as}:
\begin{gather}  
  P_C(t) = \sum_{i} \tilde{P}^C_i (t) = \sum_{i,j} P^C_{ij} (t) ,  \tilde{P}^C_i (t) = \sum_{j} P^C_{ij} (t) ,   \label{SD-diff-c} \\ 
\tilde{P}^C_i (t) = \left\{ \frac{ x_i x^T_i \mathcal{J}}{N'(\lambda_i) N(-\lambda_i)} \left(I-e^{(I \lambda_i+A^*_C) t)}\right) \right\}_H \, , \label{subG0-diff-c} \\
P^C_{ij} (t) = \left\{ \frac{e^{(\lambda_i+\lambda^*_j)t} -1}{\lambda_i+\lambda^*_j} \cdot \frac{x_i x^*_j}{N'(\lambda_i) N'(\lambda^*_j)} \right\}_H \, . 
\label{subG-diff-c} 
\end{gather}

\noindent \textbf{Remark 2.} The element-wise representation of matrices $\tilde{P}^C_i (t)$ and $P^C_{ij} (t)$ in \eqref{subG-diff-c} takes the following form:
\begin{gather} 
\tilde{P}^C_i = \sum_{\mu, \nu} \left\{\frac{(-1)^{\nu-1} \lambda_i^{\mu+\nu-2}}{N'(\lambda_i) N(-\lambda_i)} \ 
\mathbb{1}_{\mu \nu} \left(I-e^{(I \lambda_i+A^*_C) t)}\right) \right\}_H \, , \nonumber \\
P^C_{ij} = \sum_{\mu, \nu} 
\left\{\frac{e^{(\lambda_i+\lambda^*_j)t} -1}{\lambda_i+\lambda^*_j} \cdot \frac{\lambda_i^{\mu-1} (\lambda^*_j)^{\nu-1}}{N'(\lambda_i) N'(\lambda^*_j)} \mathbb{1}_{\mu \nu} \right\}_H . \nonumber  
\end{gather}\\

\noindent \textbf{Proof.} Let us find the right eigenvector $x_i=[x_1^i,x_2^i,\cdots,x_n^i]^T$ for eigenvalue $\lambda_i$ of matrix $A_C$ in \eqref{sys-c}. We choose the first component as $x_1^i=1$, and the other components are obtained from 
\begin{gather}  
A_C \, x_i= \lambda_i x_i \ \ \text{or} \ \ \nonumber \\
x^i_1 =1, \ x^i_2 = \lambda_i x^i_1 = \lambda_i, \ \nonumber \\
x^i_3 = \lambda_i x^i_2 = \lambda^2_i, \ \cdots, \ x^i_n = \lambda_i x^i_{n-1} = \lambda^{n-1}_i . \nonumber
\end{gather}  
In this case, the last equality,
\begin{equation}
- \sum^{n-1}_{k=0} a_k \lambda^k_i =\lambda_i x^i_n = \lambda^n_i \ \Leftrightarrow \ \text{det}(I\lambda_i - A_C) = 0,  \nonumber
\end{equation}
is automatically satisfied for $\lambda_i$. Thus, the right eigenvector of system \eqref{sys-c} corresponding to $\lambda_i$ can be chosen as:
\begin{equation} \label{right-ev}
x_i = [1,\lambda_i, \lambda^2_i, \cdots, \lambda^{n-1}_i]^T.
\end{equation}
The left eigenvector $y_i=[y_1^i, y_2^i,\cdots, y_n^i]^T$ corresponding to $\lambda_i$ is obtained from the condition 
\begin{gather}  
y^T_i A_C = \lambda_i y^T_i \ \ \text{or} \nonumber \\
- a_0 y^i_n = \lambda_i y^i_1 , \ y^i_1 - a_1 y^i_n = \lambda_i y^i_2 , \ \nonumber \\
y^i_2 - a_2 y^i_n = \lambda_i y^i_3 , \ \cdots, \ y^i_{n-1} - a_{n-1} y^i_n = \lambda_i y^i_n . \nonumber
\end{gather}  
We choose to normalize the component $y^i_n = -1$; then the other components of $y_i^T$ are:
\begin{gather}  
y^i_1 = \lambda^{-1}_i a_0 = - \lambda^{-1}_i \sum^n_{k=1} a_k \lambda^k_i , \nonumber \\
y^i_2 = \lambda^{-2}_i (a_0 +\lambda_i a_1) = - \lambda^{-2}_i \sum^n_{k=2} a_k \lambda^k_i , \ \cdots \nonumber \\
y^i_{n-1} = \lambda^{-(n-1)}_i \sum^{n-2}_{k=0} a_k \lambda^k_i  = - \lambda^{-(n-1)}_i (a_{n-1}\lambda^{n-1}_i +1\cdot \lambda^n_i ) , \nonumber \\
y^i_n = \lambda^{-n}_i \sum^{n-1}_{k=0} a_k \lambda^k_i  = - \lambda^{-n}_i (1\cdot \lambda^n_i ) = -1 . \nonumber
\end{gather}  
Thus, $y_i$ can be briefly written in matrix form as
\begin{equation} \label{left-ev}   
y_i =\frac{1}{\lambda^n_i} \mathcal{H}_l x_i, \ \mathcal{H}_l = 
\begin{bmatrix}
0      & \cdots    &   0     & a_0  \\ 
\vdots   & \iddots   & \iddots  & a_1 \\
 0    & a_0 & \iddots     &  \vdots \\
a_0  & a_1 &   \cdots  &a_{n-1} 
\end{bmatrix} \, ,
\end{equation}
where $\mathcal{H}_l$ is a lower Hankel matrix formed by the coefficients $a_i$ of the characteristic polynomial.
The scalar product of the right and left eigenvectors is
\begin{equation} 
y^T_i x_i = - \lambda^{-1}_i \sum^n_{k_0=1} \sum^n_{k=k_0} a_k \lambda^k_i 
= -  \sum^n_{k=1} k a_{k} \lambda^{k-1}_i = - N'(\lambda_i) . \nonumber
\end{equation}
The matrix residue of the resolvent corresponding to $\lambda_i$ can be found as
\begin{equation} \label{res-0}
R_i = \frac{x_i y_i^T}{y_i^T x_i} = \frac{x_i y_i^T}{-N'(\lambda_i)} .
\end{equation}
We find the pair sub-Gramian $P_{ij}^C$ according to \eqref{subgram-general}. Substituting into \eqref{subgram-general}
\begin{equation}
R_i = \frac{x_i y_i^T}{-N'(\lambda_i)} , \ B =e_n, \ y^T_i e_n = -1, \ e^T_n (y^T_i)^* = -1, \ \nonumber \\
\end{equation}
we obtain \eqref{subG-diff-c} after symmetrization.
Let us find the sub-Gramian $\tilde{P}^C_i$ according to \eqref{subgram-general}. In this formula
\begin{equation}
R_i B = \frac{x_i y^T_i e_n}{-N'(\lambda_i)}, \ B^T_n = e^T_n . \nonumber
\end{equation}
The vector $z_i^T=[z_1,z_2,…,z_n]=e_n^T (\lambda_i I+A_C^* )^{-1}$ we find directly from the condition
\begin{equation}
z_i^T (\lambda_i I+A_C^*) = [0,\cdots,0,1]. \nonumber
\end{equation}
Writing out this system element by element, we obtain
\begin{gather}
\lambda_i z_1+z_2=0, \ z_2 = -\lambda_i z_1 , \ \nonumber \\
\lambda_i z_2+z_3=0, \ z_3 = -\lambda_i z_2 = (-\lambda_i )^2 z_1, \ \cdots \ \nonumber \\
\lambda_i z_{n-1}+z_n=0, \ z_n = -\lambda_i z_{n-1} = (-\lambda_i )^{n-1} z_1,  \ \nonumber \\
-a_0 z_1-a_1 z_2-\cdots-a_{n-1} z_n-(-\lambda_i ) z_n = -z_1 N(-\lambda_i )=1, \ \nonumber \\
z_1 = -N^{-1} (-\lambda_i ) . \nonumber 
\end{gather}
Therefore, we obtain
\begin{equation} \label{z_i}      
z_i^T = - \frac{[1,-\lambda_i, (-\lambda_i )^2, \dots, (-\lambda_i )^{n-1}] }{N(-\lambda_i )}
=\frac{x^T_i \mathcal{J}}{N(-\lambda_i )} ,   
\end{equation}
where the matrix $\mathcal{J} = \mathrm{diag}\,(-1, (-1)^2, \dots,(-1)^n)$.
According to \eqref{subgram-general} we obtain
\begin{equation}  \label{SG-dif-c0}   
\tilde{P}^C_i = 
- \left\{\frac{x_i z^T_i}{N'(\lambda_i)} \ \left(I-e^{(I \lambda_i+A^*_C) t)}\right) \right\}_H \, , \nonumber
\end{equation}
Substituting here 
$z_i^T$ from \eqref{z_i}, we obtain \eqref{subG0-diff-c} after the symmetrization. 
The coefficients at different exponential terms in the solution of 
\eqref{Lyap-diff-c} with constant coefficients are determined uniquely. 
This implies the uniqueness of the spectral components in \eqref{subG0-diff-c} when $\lambda_i \ne \lambda_j$ 
and in \eqref{subG-diff-c} when $\lambda_i + \lambda^*_j \ne \lambda_k + \lambda^*_l$ for each $i,j,k,l = \overline{1,n}$ 
$\blacksquare$.\\

Note that it follows from the proof above that SDSE without symmetrization 
\begin{gather}
  P_C(t) = \sum_{i} \hat{P}^C_i \left(I-e^{(I \lambda_i+A^*_C) t)}\right)  , \nonumber \\
  \hat{P}^C_i = \frac{x_i x^T_i \mathcal{J}}{N'(\lambda_i) N(-\lambda_i)} \label{subG-diff-c-non-sym}
\end{gather}
is also valid, and in some cases, it may be more convenient.\\

\noindent \textbf{Remark 3.} 
If a non-zero boundary condition $P_C (0)=P_C^0$ is specified in \eqref{Lyap-diff-c}, then the general solution to \eqref{Lyap-diff-c} consists of two parts
\begin{equation}
\tilde{P}_C (t) = P_C (t) + P_C^0 (t) , \nonumber
\end{equation}
where the expansion of $P_C (t)$ is defined by Theorem 2, and $P_C^0 (t)$ satisfies Lyapunov equation \eqref{Lyap-diff-c} with a zero right-hand side and the initial condition $P_C^0 (0) = P_C^0$. In this case, the spectral decompositions of the second part of the solution have a similar form:
\begin{align}                        
P_C^0 (t) = & \sum_{i} \tilde{P}^C_{0\,i} (t) \nonumber \\
& = \sum_{i} \left\{\frac{x_i y^T_i}{-N'(\lambda_i)} P^0_C \ e^{(I \lambda_i+A^*_C) t} \right\}_H \, , \label{SD-inicon-1} \\
P_C^0 (t) = & \sum_{i} P^C_{0\,ij} (t) \nonumber \\
& = \sum_{i,j} \left\{\frac{x_i y^T_i P^0_C (y^T_j)^* x^*_j }{N'(\lambda_i) N'(\lambda^*_j)}
 \ e^{(\lambda_i+\lambda^*_j) t} \right\}_H \, . \label{SD-inicon-2} 
\end{align}

\noindent \textbf{Proof.} Using \eqref{res-0} and the following properties of residues 
\begin{equation}
A_C R_i = \lambda_i R_i, \ R_j^* A_C^* = R_j^* \lambda_j^*, \ \sum_i R_i = \sum_j R^*_j = I , \nonumber 
\end{equation}
we can directly check the validity of decomposition \eqref{SD-inicon-1}:
\begin{gather}
\text{for each }i = \overline{1,n}: \ \  - \frac{d\tilde{P}^C_{0\,i}(t)}{dt} +A_C \, \tilde{P}^C_{0\,i}(t) + \tilde{P}^C_{0\,i}(t) A^T_C = 0  \ \nonumber \\ 
\Rightarrow \ \  - \frac{d P^C_0(t)}{dt} +A_C \, P^C_0(t) + P^C_0(t) A^T_C = 0 , \ \nonumber \\ 
P^0_C(0) = \sum_{i} \tilde{P}^C_{0\,i} (0) = \sum_i R_i P^0_C = P^0_C . \nonumber 
\end{gather}
Similarly, we check the validity of \eqref{SD-inicon-2} $\blacksquare$.\\

\noindent \textbf{Example 2} illustrates SDSE of the finite Gramian in companion form. 
According to Theorem 2 for the system considered in Example 1, 
the solution \eqref{Lyap-diff-c} with zero boundary condition $P_C (0)=0$ 
is decomposed over the pair spectrum in the form
\begin{gather}
P_C (t) = P_{11}^C (1-e^{2t}) + (P_{12}^C+P_{21}^C )(1-e^{3t} ) \nonumber \\
+ \ (P_{22}^C+P_{31}^C+P_{13}^C )(1-e^{4t} )  \nonumber \\
+ \ (P_{23}^C+P_{32}^C )(1-e^{5t} )+P_{33}^C (1-e^{6t} ) ,  \nonumber   
\end{gather}
where the matrices $P_{ij}^C$ were found in Example 1. 
If the initial condition, $P_C (0)=P_C^0$ is given in \eqref{Lyap-diff-c}, 
then according to \eqref{SD-inicon-2}, an additional term $P_C^0(t)$ is added to the solution 
\begin{gather}
\tilde{P}_C (t) = P_C (t) + P_C^0 (t), \ \nonumber \\
P_C^0 (t) = \sum^3_{i,j=1} e^{(\lambda_i+\lambda_j^* ) t}  \left\{R_i P_C^0 R_j^* \right\}_H , \nonumber
\end{gather}
where the eigenvalues and residue matrices \eqref{res-0} are
\begin{gather}
\lambda_1=1, \ \lambda_2=2, \ \lambda_3=3, \  R_1=\frac{1}{2}
\begin{bmatrix}
6 & -5 & 1 \\
6 & -5 & 1 \\
6 & -5 & 1    
\end{bmatrix} , \nonumber \\
R_2 =
\begin{bmatrix}
-3 & 4 & -1 \\
-6 & 8 & -2 \\
-12 & 16 & -4    
\end{bmatrix} , \ R_3 = \frac{1}{2}
\begin{bmatrix}
2 & -3 & 1 \\
6 & -9 & 3 \\
18 & -27 & 9    
\end{bmatrix} \ \square .  \nonumber 
\end{gather}\\

\section{Relationship between SDSE of Gramian in Companion Form and 
Gramian of MI System}

In this section, we show in Theorem 3 how to obtain the SDSE of the controllability Gramian for an arbitrary controllable MI system \eqref{sys-0} 
from the corresponding SDSE previously obtained in the canonical controllability form. To find the spectral components of the controllability Gramian in \eqref{Lyap-diff} for the SI controllable system \eqref{sys-0}, it is necessary to use the matrix $T$ of the transition to the original coordinate system \eqref{transfer-matrix} 
by analogy with expression \eqref{P-to-orig-sys}
\begin{gather}
P=\sum_{i} \tilde{P}_i = \sum_{i,j} P_{ij}, \ \nonumber \\
\tilde{P}_i = T \tilde{P}_i^C T ^T, \ \  P_{ij} = T P_{ij}^C \,T^T  . \nonumber
\end{gather}
In the case of an MI system with $m$ inputs, the input matrix $B$ can be split into separate columns $b_\gamma$. Then the right-hand side of \eqref{Lyap-diff} is partitioned into parts corresponding to individual SI systems, and the controllability Gramian is partitioned into independent parts corresponding to them:
\begin{gather}
AP+P A^T= -BB^T,  \ \ BB^T=\sum_{\gamma=1}^m b_{\gamma} b_{\gamma}^T  \ \nonumber \\
\Rightarrow \ \ P=\sum^m_{\gamma=1} P_{\gamma}, \   AP_{\gamma} + P_{\gamma} A^T = -b_{\gamma} b_{\gamma}^T . \nonumber
\end{gather}
By analogy with \eqref{MI-P-from-Pc} for the controllability Gramian, we can obtain similar expressions for its spectral components.\\

\noindent \textbf{Theorem 3.} \textit{Let MI system \eqref{sys-0} have a simple spectrum. Then SDSE of the finite controllability Gramian $P(t)$, satisfying \eqref{Lyap-diff} with zero initial condition $P(0)=0$, can be expressed through the sub-Gramians $\tilde{P}_i^C (t)$ and $P_{ij}^C (t)$ obtained in the corresponding canonical system as}
\begin{gather}
P(t) = \sum_{i} \tilde{P}_i (t) = \sum_{i,j} P_{ij} (t) , \ \nonumber \\ 
\tilde{P}_i (t) = \mathcal{C} (\mathcal{H}_u \tilde{P}_i^C (t) \mathcal{H}_u \otimes I_m ) \mathcal{C}^*,  \  \label{subgram-diff-MI}  \\   
P_{ij} (t) = \mathcal{C} (\mathcal{H}_u P_{ij}^C (t) \mathcal{H}_u \otimes I_m ) \mathcal{C}^*, \nonumber
\end{gather}
\textit{where $\mathcal{C} = [B,AB,A^2 B,\dots,A^{n-1} B]$ is the controllability matrix, $I_m$ is a unit matrix $m\times m$, $\mathcal{H}_u$ 
is given by \eqref{transfer-matrix}, and the sub-Gramians $\tilde{P}_i^C (t)$ and $P_{ij}^C (t)$ 
are obtained by \eqref{subG-diff-c} in Theorem 2.} \\ 

\noindent \textbf{The proof} follows directly from Theorem 2 and formula \eqref{MI-P-from-Pc} derived for the Gramian in \cite{Hauksdottir-2009}. 
The computation of $\tilde{P}_i^C (t)$ and $P_{ij}^C (t)$ in \eqref{subgram-diff-MI} depends only on the dynamics matrix $A$ 
and is independent of the input matrix $B$.\\  

\noindent \textbf{Remark 4.} If a non-zero initial condition $P(0)=P_0$ is specified in \eqref{Lyap-diff}, then its solution consists of two parts
\begin{equation}
\tilde{P}(t) = P(t) + P_0 (t) , \nonumber
\end{equation}
where the decomposition of $P(t)$ is determined by Theorem~3, and $P_0 (t)$ satisfies \eqref{Lyap-diff} with a zero right-hand side and initial condition $P(0)=P_0$. 
In this case, the spectral decomposition and pairwise spectral decomposition of $P_0 (t)$ are computed similarly to \eqref{SD-inicon-1} and \eqref{SD-inicon-2} as
\begin{gather}
P_0 (t) = \sum_{i} \tilde{P}_{0\, i} (t) = \sum_{i,j} P_{0\, ij} (t)  , \ \nonumber \\   
\tilde{P}_{0\, i} (t) = \left\{R_i P_0 e^{(I \lambda_i + A_C^* )t} \right\}_H , \ \nonumber \\ 
P_{0\, ij} (t) = \left\{R_i P_0 R_j^* e^{(\lambda_i+\lambda_j^* )t} \right\}_H,  \nonumber      
\end{gather}
where $R_i, R_j$ are the residue matrices, which can be calculated more easily using eigenvectors in the original coordinate system.\\

\section{Spectral Decomposition of Inverse of Gramian in Companion Form}

In this section, we derive the SDSE of the inverse of the controllability Gramian in Theorem~4 and show in Property~4 that the corresponding spectral components in this decomposition are uniquely defined as measurable quantities in the minimum energy control problem. 
As an additional result in Corollary~1, we obtain a closed-form solution for a certain class of Riccati equations.
In Theorem 5, we derive the SDSE of the inverse of the {\it finite} controllability Gramian, which is a solution to the differential Riccati equation, 
and summarize some results of Theorems 1, 2, and 4.

Let us consider system \eqref{sys-c} in the controllability canonical form. 
If the infinite Gramian $P_{C}$ satisfies the algebraic Lyapunov equation \eqref{Lyap-algeb}, then its inverse 
$P_{C}^{-1}$ satisfies the algebraic Riccati equation
 \begin{equation}
P_C^{-1} A_C + A_C ^T P_C^{-1} = - P_C^{-1} b_C \, (b_C )^T P_C^{-1} .  \label{Ric-algeb-c}
\end{equation}
The matrix $P_C^{-1}$ can be  
represented as a sum of Hermitian matrices corresponding to individual eigenvalues
or their pairwise combinations. \\ 

\noindent \textbf{Theorem 4.} \textit{Let system 
\eqref{sys-c} have a simple spectrum. 
Then SDSE of the inverse of the \textbf{infinite} Gramian $P_C^{-1}$ in \eqref{Ric-algeb-c} is given by}
\begin{gather}    
P_C^{-1} = \sum_{j} \tilde{P}_j^{-C} = \sum_{j} \hat{P}_j^{-C} = \sum_{i,j} \hat{P}_{ij}^{-C}, 
\ \hat{P}_j^{-C} = \sum_{i} \hat{P}_{ij}^{-C},\ \nonumber \\
\tilde{P}_j^{-C} = \left\{ \hat{P}_j^{-C} \right\}_H, \ P_{ij}^{-C} = \left\{ \hat{P}_{ij}^{-C} \right\}_H, \ \nonumber \\
\hat{P}_j^{-C} = \frac{N(-\lambda_j)}{-N'(\lambda_j)} \mathcal{J}  y_j  y_j^T , \ \label{SD-invG} \\
\hat{P}_{ij}^{-C} = \frac{N(-\lambda^*_i) N(-\lambda_j)}{-N'(\lambda^*_i) N'(\lambda_j)} \cdot \frac{(y^*_i)^T y^T_j}{\lambda^*_i+\lambda_j}, \label{SD-invG-2} \\
 y_j=\frac{\mathcal{H}_l  x_j}{\lambda_j^n} , \ 
 x_j = [1,\lambda_j, \lambda^2_j, \cdots, \lambda^{n-1}_j ]^T \ , \nonumber
\end{gather}
\textit {where $\mathcal{H}_l$ is defined in \eqref{left-ev}, 
and matrices $\hat{P}_j^{-C}$ and $\hat{P}_{ij}^{-C}$ in \eqref{SD-invG}, \eqref{SD-invG-2} are uniquely defined by the condition} 
\begin{equation}
\forall \ i ,j: \ \ 
\hat{P}_i^{C} \hat{P}_j^{-C} = \delta_{ij} R_i , \  \hat{P}_i^{C} = \frac{x_i x^T_i \mathcal{J}}{N'(\lambda_i) N(-\lambda_i)} \label{orto-cond} 
\end{equation}
\textit {of orthogonality of eigenparts in the expansions of $P_C$ and $P_C^{-1}$ in \eqref{subG-alg-c} and \eqref{SD-invG}.}\\

\noindent \textbf{Remark 5.} The element-wise representation of matrices $\hat{P}_j^{-C}$ and $\hat{P}_{ij}^{-C}$ in \eqref{SD-invG} and \eqref{SD-invG-2} takes the following form:
\begin{gather}  
(\hat{P}_j^{-C})_{\mu \nu} = \frac{(-1)^{\mu-1} N(-\lambda_j)}{N'(\lambda_j) \lambda_j^{\mu+\nu} } \ \sum_{k=\mu}^n a_k \lambda_j^k \
 \sum_{l=\nu}^n a_l \lambda_j^l , \nonumber \\
(\hat{P}_{ij}^{-C})_{\mu \nu} = \frac{-N(-\lambda^*_i) N(-\lambda_j)}
{N'(\lambda^*_i) N'(\lambda_j)(\lambda^*_i+\lambda_j)(\lambda^*_i)^{\mu}\lambda^{\nu}_j} \cdot \nonumber \\
\cdot \sum_{k=\mu}^n a_k (\lambda^*_i)^k  \ \sum_{l=\nu}^n a_l \lambda_j^l .  \label{SD-invG-elem}
\end{gather}

\noindent \textbf{Proof.} 
For the time-independent part of the controllability Gramian $\hat{P}_C$ in \eqref{subG-diff-c-non-sym} 
\begin{equation}
P_C = \sum_{i} \hat{P}_i^{C} = \sum_{i} \frac{x_i x^T_i \mathcal{J}}{N'(\lambda_i) N(-\lambda_i)} \nonumber
\end{equation}
we construct the inverse matrix in the form 
\begin{equation} 
P_C^{-1}= \sum_{j} \hat{P}_j^{-C} = \sum_j \frac{ \mathcal{J} y_j y_j^T N(-\lambda_j)}{-N'(\lambda_j )} .     \label{invG-0}   
\end{equation}
Using the relations $ \mathcal{J}^2 = I$, $x^T_i y_i = -N'(\lambda_i)$ and orthogonality of eigenvectors $x_i$ and $y_j$ when $i\ne j$, we obtain
\begin{gather}
P_C P_C^{-1} = \sum_{i,j} \hat{P}_i^C \hat{P}_j^{-C}  \nonumber \\
= \sum_{i,j} \frac{x_i x_i^T \mathcal{J}}{N'(\lambda_i)N(-\lambda_i)} \cdot \frac{ \mathcal{J} y_j y_j^T N(-\lambda_j )}{N'(\lambda_j)} \nonumber \\
= \sum_{i,j} \frac{x_i \delta_{ij} y_j^T N(-\lambda_j)}{-N'(\lambda_i)N(-\lambda_i)} = \sum_i R_i = I . \nonumber
\end{gather}
Carrying out the symmetrization
\begin{equation}
P_C^{-1}=\{P_C^{-1} \}_H = \sum_{j} \{\hat{P}_j^{-C} \}_H = \sum_j \tilde{P}_j^{-C} ,  \nonumber
\end{equation}
we obtain \eqref{SD-invG}. Let us prove that 
\begin{equation}
\hat{P}_{ij}^{-C} = R^*_i \hat{P}_{j}^{-C} . \label{th4-1} 
\end{equation}
Substituting here \eqref{SD-invG}, \eqref{SD-invG-2} and $R^*_i = \frac{(y^*_i)^T x^*_i}{-N'(\lambda^*_i)}$ we obtain that \eqref{th4-1} holds if
\begin{equation}
\frac{-N(-\lambda^*_i)}{\lambda^*_i+\lambda_j} = x^*_i \mathcal{J} y_j . \label{th4-2} 
\end{equation}
On the other hand, from \eqref{z_i} we have
\begin{gather}
\forall \lambda_i \in \sigma(A_C) : \ \ \frac{x^T_i \mathcal{J}}{N(-\lambda_i)} = e^T_n (I \lambda_i+A^*_C)^{-1} \nonumber \\
\Rightarrow \ \ \lambda^*_i \in \sigma(A_C) : \ \ \frac{x^*_i \mathcal{J}}{N(-\lambda^*_i)} = e^T_n (I \lambda^*_i +A^*_C)^{-1} . \label{th4-3} 
\end{gather}
Multiplying \eqref{th4-3} on the right by $y_j$, and taking into account that from \eqref{left-ev} it follows that
\begin{gather}
 e^T_n y_j = -1, \ (I \lambda^*_i +A^*_C) y_j = (\lambda^*_i+\lambda_j) y_j, \nonumber \\
\Rightarrow \ \ e^T_n (I \lambda^*_i +A^*_C)^{-1} y_j = \frac{e^T_n y_j }{\lambda^*_i+\lambda_j} = \frac{-1}{\lambda^*_i+\lambda_j} , \nonumber
\end{gather}
we thus obtain \eqref{th4-2}. Therefore \eqref{th4-1} is verified. From \eqref{th4-1}, we obtain that $\hat{P}_j^{-C} = \sum_{i} \hat{P}_{ij}^{-C}$ and the validity of \eqref{SD-invG-2}. Substituting into \eqref{SD-invG} and \eqref{SD-invG-2} the expressions \eqref{left-ev} for the components of $y_j$, we obtain \eqref{SD-invG-elem}. For matrices $\hat{P}_j^{-C}$ in \eqref{SD-invG} the condition \eqref{orto-cond} clearly holds. 
Let us prove the uniqueness of the representation of \eqref{SD-invG}. 
It follows from \eqref{orto-cond} that $\hat{P}_j^{-C}$ is orthogonal to $(n-1)$ independent vectors $x^T_i \mathcal{J}, \ i\ne j$, 
that is, $\mathrm{rank} \ \hat{P}_j^{-C} = 1$. So the following representation $ \hat{P}_j^{-C} = a \, b^T $
through some non-zero vectors $a, b \in \mathbb{C}^n$ is valid.
The vector $a$ must be orthogonal to $(n-1)$ vectors $x^T_i \mathcal{J}, \ i\ne j$. 
Hence, there must be $a = \alpha \cdot \mathcal{J} y_j$ for some non-zero $\alpha \in \mathbb{C}$. There must also be 
\[
 \hat{P}_j^{C} \hat{P}_j^{-C} =  \frac{x_j x^T_j \mathcal{J}}{N'(\lambda_j) N(-\lambda_j)} \cdot \alpha \mathcal{J} y_j b^T 
=\frac{- \alpha x_j b^T} {N(-\lambda_j)} = R_j \, .
\]
Substituting here \eqref{res-0} for $R_j$, we obtain 
\[
b^T =  \frac{y^T_j N(-\lambda_j)}{-\alpha N'(\lambda_j)} ,
\]
and $\hat{P}_j^{-C}$ must satisfy \eqref{SD-invG}. 
The uniqueness of \eqref{SD-invG} is proved. Then the uniqueness of \eqref{SD-invG-2} follows from \eqref{th4-1} $\blacksquare$.\\

The spectral components of $\tilde{P}^{-C}_j$ and $P^{-C}_{ij}$ in Theorem~4 are not merely mathematical entities but also measurable quantities. 
They can be interpreted using the problem \eqref{Emin} of {\it minimum energy control}.\\

\noindent \textbf{Property 4 (Interpretation).}
\textit{Let system \eqref{sys-c} be stable and state $x_0$ be reachable.
Then the eigenparts $\tilde{P}^{-C}_i$ and $P^{-C}_{ij}$ in \eqref{SD-invG} and \eqref{SD-invG-2} 
characterize, respectively, the linear and quadratic partitions of the minimum control energy $E_{min}$ in \eqref{Emin} 
into real parts corresponding to the components of $x_0$ in the eigenbasis of matrix $A_C$, i.e.,}
\begin{gather}
E_{min} = x^T_0 P^{-1}_C x_0 = \sum_i E_i =\sum_{ij} \hat{E}_{ij},  \ \ \text{where} \nonumber \\
E_i = x^T_0 \tilde{P}^{-C}_i x_0 = \frac{1}{2} \left( x^*_{i0} P^{-1}_C x_0 + x^*_0 P^{-1}_C x_{i0} \right) , \ \label{SD-Emin} \\
\hat{E}_{ij} = x^T_0 {P}^{-C}_{ij} x_0 = \frac{1}{2} \left( x^*_{i0} P^{-1}_C x_{j0} + x^*_{j0} P^{-1}_C x_{i0} \right), \ \label{SD-Emin-2} \\
x_0 = \sum_i x_{i0}, \ \ x_{i0} = R_i x_0 ,  \ \ x_{j0} = R_j x_0 , \nonumber 
\end{gather}
\textit{and $R_i, R_j$ are the matrix residues given by \eqref{res-0}.} \\

\noindent \textbf{Proof.} Using the orthogonality of eigenvectors, we obtain
\begin{gather}
P^{-1}_C R_i = \sum_j \frac{\mathcal{J} y_j y_j^T N(-\lambda_j)}{-N'(\lambda_j )} \cdot \frac{x_i y_i^T}{-N'(\lambda_i)} \nonumber \\
= \frac{\mathcal{J} y_i y_i^T N(-\lambda_i)}{-N'(\lambda_i)} \cdot \frac{x_i y_i^T}{-N'(\lambda_i)} = \hat{P}^{-C}_i ; \nonumber \\
\hat{P}^{-C}_i x_0 = \sum_j \hat{P}^{-C}_i R_j x_0 = \sum_j \delta_{ij} \hat{P}^{-C}_i R_j \, x_0 \nonumber \\
= \hat{P}^{-C}_i R_i \, x_0 =  \hat{P}^{-C}_i x_{i0} \ ; \nonumber \\
\hat{P}^{-C}_i x_0 = \hat{P}^{-C}_i x_{i0} = P^{-1}_C R_i \, x_{i0} = P^{-1}_C x_{i0} . \nonumber
\end{gather}
Using these relations and $\sum_i R_i = \sum_i R^*_i  = I$, we get
\begin{gather}
E_{min} = x^T_0 P^{-1}_C x_0 = \frac{1}{2} \sum_i \left( x^T_0 R^*_i P^{-1}_C x_0 + x^T_0 P^{-1}_C  R_i \,x_0\right)  \nonumber \\
=  \frac{1}{2} \sum_i \left( x^T_0 (\hat{P}^{-C}_i)^* x_0 + x^T_0 \hat{P}^{-C}_i x_0\right) \nonumber \\
= \sum_i x^T_0 \tilde{P}^{-C}_i x_0 = \frac{1}{2} \sum_i \left( x^*_{i0} P^{-1}_C x_0 + x^T_0 P^{-1}_C x_{i0} \right) .  \nonumber
\end{gather}
We similarly prove \eqref{SD-Emin-2}. $\blacksquare$\\

This interpretation of eigenparts of the inverse of the controllability Gramian carries over to an arbitrary controllable system \eqref{sys-0} 
after transition to the corresponding coordinate system. Therefore, the symmetrized eigenparts $\tilde{P}^{-C}_j$ in \eqref{SD-invG} and $P^{-C}_{ij}$ in \eqref{SD-invG-2} are uniquely determined as measurable quantities in the minimum energy control problem \eqref{Emin}. 

Similarly to the proof of Properties 2 and 3, the following properties are proved from \eqref{SD-invG-elem} and \eqref{SD-invG-2}.\\

\noindent \textbf{Property 5 (Structure).}  \textit{The inverse of the controllability Gramian $P_C^{-1}$ in the controllability canonical form 
and its spectral components $\tilde{P}_j^{-C}$ in \eqref{SD-invG} have zero-plaid structure \eqref{zero-plaid}, 
i.e., $ \left(\tilde{P}_j^{-C} \right)_{\mu \nu} = 0$ when $\mu + \nu = 2m+1$.}\\

\noindent \textbf{Property 6.} \textit{If system \eqref{sys-c} is stable, then 
$P^{-C}_{ii} \ge 0$, $i = \overline{1,n}$.}\\

\noindent \textbf{Remark 6.} It also follows from \eqref{SD-invG-elem} that computing the components of any matrix $\tilde{P}_j^{-C}$ 
or $P^{-C}_{ij}$ requires at most $O(n^2)$ operations. Indeed, for each $\lambda_i$, the components of the eigenvectors $x_i$ and $y_i$ 
are computed according to \eqref{right-ev} and \eqref{left-ev} recursively in $O(n)$ operations 
and immediately substituted into \eqref{SD-invG}, \eqref{SD-invG-2} or \eqref{SD-invG-elem}.\\

\noindent \textbf{Remark 7.} In Theorems 2 and 4, both symmetrized terms and unsymmetrized terms can be used in the spectral expansions.
However, it seems that only the symmetric parts of the spectral components have a physical interpretation, and the non-symmetric parts can be chosen conveniently as long as their sum is zero. In particular, if it is desirable to preserve the orthogonality of the spectral components, then one can 
use $\hat{P}_j^{-C}$ and $\hat{P}_{ij}^{-C}$ in \eqref{SD-invG} and \eqref{SD-invG-2},
which are not symmetric but are uniquely defined. \\ 

\noindent \textbf{Example 3} illustrates SDSE of the inverse of Gramian. 
Consider the system from Example 1. The inverse of its controllability Gramian can be calculated in closed form as
 \begin{equation}
P_C^{-1} = -12 \cdot
\begin{bmatrix}
11 & 0 & 1 \\
0 & 10 & 0 \\
1 & 0 & 1    
\end{bmatrix} = \tilde{P}_1^{-C} + \tilde{P}_2^{-C} + \tilde{P}_3^{-C},  \nonumber
\end{equation}
where the spectral components are found by \eqref{SD-invG} or \eqref{SD-invG-elem}
\begin{gather}
\tilde{P}_1^{-C} = 12
\begin{bmatrix}
-36 & 0 & -6 \\
  0 & 25 & 0 \\
 -6 & 0  & -1    
\end{bmatrix} , \nonumber \\
\tilde{P}_2^{-C} = 60
\begin{bmatrix}
9 & 0 & 3 \\
0 & -16 & 0 \\
3 &   0 & 1    
\end{bmatrix} , \ \tilde{P}_3^{-C} = 60
\begin{bmatrix}
-4 & 0 & -2 \\
 0 &  9 & 0 \\
-2 &  0 & -1    
\end{bmatrix}  . \nonumber 
\end{gather} 
The pair spectral components are obtained by \eqref{SD-invG-2}
\begin{gather}
P^{-C}_{11} = -72
\begin{bmatrix}
36 & -30 & 6 \\
-30 & 25 & -5 \\
6 & -5 & 1    
\end{bmatrix} , \ \nonumber \\
P^{-C}_{12} = P^{-C}_{21} = 120
\begin{bmatrix}
36 & -39 & 9 \\
-39 & 40 & -9 \\
9 & -9 & 2    
\end{bmatrix} , \ \nonumber \\
P^{-C}_{22} = -900
\begin{bmatrix}
9 & -12 & 3 \\
-12 & 16 & -4 \\
3 & -4 & 1   
\end{bmatrix} , \ \nonumber \\
P^{-C}_{23} = P^{-C}_{32} = 360
\begin{bmatrix}
12 & -17 & 5 \\
-17 & 24 & -7 \\
5 & -7 & 2    
\end{bmatrix} , \nonumber \\
P^{-C}_{33} = -600
\begin{bmatrix}
4 & -6 & 2 \\
-6 & 9 & -3 \\
2 & -3 & 1    
\end{bmatrix} , \ \nonumber 
\end{gather}
\begin{gather}
P^{-C}_{31} = P^{-C}_{13} = -180
\begin{bmatrix}
12 & -14 & 4 \\
-14 & 15 & -4 \\
4 & -4 & 1    
\end{bmatrix} . \nonumber 
\end{gather}
One can verify that 
\begin{gather}
\tilde{P}^{-C}_1 = P^{-C}_{11} + P^{-C}_{12} + P^{-C}_{13},  \ \tilde{P}^{-C}_2 = P^{-C}_{21} + P^{-C}_{22} + P^{-C}_{23}, \nonumber \\
\tilde{P}^{-C}_3 = P^{-C}_{31} + P^{-C}_{32} + P^{-C}_{33},  \ \ P^{-1}_C = \sum^n_{i,j=1} P^{-C}_{ij} \ \square \ . \nonumber 
\end{gather}

Consider the inverse of the controllability Gramian of an arbitrary controllable SI system.
Using Theorem 4 and transformations \eqref{P-to-orig-sys} and \eqref{transfer-matrix} for the transition to the original basis, 
we can immediately obtain the following corollary.\\

\noindent \textbf{Corollary 1.} \textit{The solution of the Riccati equation}
\begin{equation}
P^{-1} A + A^* P^{-1} = - P^{-1} b\,b^T P^{-1}  \label{Riccati} 
\end{equation}
\textit{for matrix $A$ with a simple spectrum $\sigma(A)=\{\lambda_1,\lambda_2,\cdots,\lambda_n\}$ 
such that $\lambda_i+\lambda_j \ne 0$ for all $\lambda_i,\lambda_j \in \sigma(A)$ and non-singular matrix
$\mathcal{C} = [ b,Ab,\dots,A^{n-1} b]$, has the following spectral decomposition into Hermitian matrices corresponding 
to individual eigenvalues or their pairwise combinations}
\begin{gather}
P^{-1} = \sum_{j} (\mathcal{C}^T)^{-1} \mathcal{H}_u ^{-1} \tilde{P}_j^{-C} \mathcal{H}_u^{-1} \mathcal{C}^{-1} , \ \nonumber \\
P^{-1} = \sum_{i,j} (\mathcal{C}^T)^{-1} \mathcal{H}_u ^{-1} P_{ij}^{-C} \mathcal{H}_u^{-1} \mathcal{C}^{-1} , \label{Riccati-sol} 
\end{gather}
\textit{where $\mathcal{H}_u$ is given by \eqref{transfer-matrix}, and matrices $\tilde{P}_j^{-C}$ and $P^{-C}_{ij}$ are defined by \eqref{SD-invG} and \eqref{SD-invG-2} in Theorem $4$.}\\

If the finite Gramian $P_C(t)$ satisfies the differential equation \eqref{Lyap-diff-c}, then its inverse $P_C^{-1}(t)$ satisfies the differential Riccati equation
\begin{equation}
\frac{dP_C^{-1}}{dt} + P_C^{-1} A_C + A_C ^T P_C^{-1} = - P_C^{-1} b_C \, (b_C )^T P_C^{-1} .  \label{Ric-diff-c}                  
\end{equation}

To study the spectral properties of the inverse of Gramian on an arbitrary time interval and to analyze the dynamics of the system at future times, we derive the SDSE of the inverse of the finite controllability Gramian, which is a solution to the differential Riccati equation \eqref{Ric-diff-c}. In Theorem 5, we also summarize some previous results of Theorems 1, 2, and 4. \\

\noindent \textbf{Theorem 5.} \textit{Let system \eqref{sys-c} 
have a simple spectrum.
Then SDSE of \textbf{the finite} controllability Gramian $P_C (t)$ and its inverse $P_C^{-1} (t)$, 
satisfying \eqref{Lyap-diff-c}, \eqref{Ric-diff-c}, and the boundary condition $P_C (0)=P_0$, takes the following form}:
\begin{gather}
P_C (t) = P_{C}(\infty) +\sum_i \left(R_i P_0 - \hat{P}_i^C\right) \, e^{(\lambda_i I+A_C^T ) t} , \nonumber \\ 
P_C^{-1}(t) = G(t) \sum_j \hat{P}_j^{-C}, \ \label{SD-invFG} \\
\text{where} \ \ \ P_{C}(\infty) = \sum_i \hat{P}_i^C, \ \ P^{-1}_{C}(\infty) = \sum_j \hat{P}_j^{-C} , \nonumber \\ 
\hat{P}_i^C = \frac{x_i x_i^T \mathcal{J}}{-N'(\lambda_i) N(-\lambda_i)}, \ 
\hat{P}_j^{-C} = \frac{\mathcal{J} y_j y_j^T N(-\lambda_j )}{-N'(\lambda_j)}, \ \label{inv-SG}  \\ 
 G^{-1}(t) = I - e^{(A_C^T +\mathcal{J} A_C^T \mathcal{J}) t} + \sum_i \hat{P}_i^{-C}P_0 e^{(\lambda_i I+A_C^T ) t}, \nonumber
\end{gather}
\textit{$x_i,y_j,R_i$ are defined in \eqref{right-ev}, \eqref{left-ev}, \eqref{res-0}, 
$G(t)$ is the normalization matrix, and $\mathcal{J}=\mathrm{diag} \, (-1,(-1)^2,\dots,(-1)^{n-1})$.}\\

\noindent \textbf{Proof.} SDSE for the solutions of algebraic Lyapunov and Riccati equations $P_{C}(\infty)$ and $P^{-1}_C(\infty)$ were obtained
in Theorems 1 and 4. It was shown in Theorem 2 and Remark 3 that the first expansion for $P_C (t)$ satisfies \eqref{Lyap-diff-c}. 
At $t=0$ one can directly check that $P_C (0) = \sum_i R_i P_0 = P_0$.
Given that $\frac{N(-\lambda_j ) y_j^T x_i}{-N'(\lambda_i) N(-\lambda_i )} = \delta_{ij}$ is the Kronecker symbol, 
we observe that 
\begin{gather}
\hat{P}_j^{-C} \hat{P}_i^C = \frac{\mathcal{J} y_j y_j^T N(-\lambda_j )}{-N'(\lambda_j)}
\cdot \frac{x_i x_i^T \mathcal{J}}{-N'(\lambda_i) N(-\lambda_i)} \nonumber \\
= \frac{\mathcal{J} y_j \delta_{ji} x_i^T \mathcal{J}}{-N'(\lambda_j)} = \mathcal{J} R_j^T \delta_{ji} \mathcal{J}, \ \nonumber \\
\hat{P}_j^{-C} R_i = \frac{\mathcal{J} y_j y_j^T N(-\lambda_j )}{-N'(\lambda_j)} \cdot \frac{x_i y_i^T}{-N'(\lambda_i)} \ \nonumber \\
= \frac{\mathcal{J} y_j N(-\lambda_j ) \delta_{ji} y^T_j}{-N'(\lambda_j)} = \delta_{ji} \hat{P}_j^{-C} \ . \nonumber
\end{gather}
We check the expansion for $P_C^{-1}(t)$ by direct substitution:
\begin{gather}
G^{-1}(t) P_C^{-1} (t) P_C (t) = \sum_{j,i} \hat{P}_j^{-C} \hat{P}_i^C (I - e^{(\lambda_i I+A_C^T ) t}) \nonumber \\
+ \ \sum_{j,i} \hat{P}_j^{-C} R_i P_0 e^{(\lambda_i I+A_C^T ) t} \nonumber \\
=  \sum_{j,i} \mathcal{J} R_j^T \delta_{ji} \mathcal{J} (I - e^{(\lambda_i I+A_C^T ) t})
+ \delta_{ji} \hat{P}_j^{-C} P_0 e^{(\lambda_i I+A_C^T ) t} \nonumber \\
= \mathcal{J}  \sum_{i} R_i^T \mathcal{J} - \sum_i \mathcal{J} R^T_i \mathcal{J} e^{\lambda_i t} e^{A_C^T t} \nonumber \\
+ \sum_i \hat{P}_i^{-C}P_0 e^{(\lambda_i I+A_C^T ) t} 
= G^{-1}(t) \ \ \ \ \blacksquare . \nonumber
\end{gather}

\noindent \textbf{Remark 8.} Since $P_C=\{P_C\}_H$, $P_C^{-1}=\{P_C^{-1} \}_H$, then the decomposition \eqref{SD-invFG} 
at $P_0=0$ can be symmetrized as follows:
\begin{gather}
P_C (t) = \sum_i \left\{ \frac{-x_i x_i^T \mathcal{J}}{N'(\lambda_i) N(-\lambda_i)} (I - e^{(\lambda_i I+A_C^T ) t})\right\}_H , \ \nonumber \\
P_C^{-1} (t) = \sum_j \left\{ G(t) \frac{\mathcal{J} y_j y_j^T N(-\lambda_j)}{-N'(\lambda_j)} \right\}_H.  \nonumber 
\end{gather}
However, for the symmetrized terms of the spectral expansion, the property of their orthogonality is no longer satisfied.\\

\noindent \textbf{Remark 9.} 
Quadratic SDSEs of the Gramian and its inverse under the conditions of Theorem 5 are obtained similarly as
\begin{gather}
P_C (t) = \sum_{i,j} P_{ij}^C +\sum_{i,j} \left(R_i P_0 R^*_j - P_{ij}^C\right) \, e^{(\lambda_i + \lambda^*_j ) t} , \nonumber \\ 
P_C^{-1}(t) = G(t) \sum_{ij} P_{ij}^{-C}. \ \label{SD-invFG-2} 
\end{gather}

\noindent \textbf{Example 4} illustrates an application of Theorem 5. Consider the system from Example 1. 
Let us denote $F_i(t) = I - e^{(\lambda_i I+A_C^T ) t}$.
According to \eqref{right-ev}, \eqref{left-ev}, and \eqref{res-0}:
\begin{gather}
x_1 = [1,1,1]^T, \ x_2=[1,2,4]^T, \ x_3=[1,3,9]^T , \nonumber \\ 
y_1 = [-6,5,-1]^T, \ y_2=[-3,4,-1]^T, \ y_3 = [-2,3,-1]^T , \nonumber \\
N'(\lambda_1) = 2, \ N'(\lambda_2) = -1, \ N'(\lambda_3) = 2 . \nonumber
\end{gather}
The unsymmetrized spectral components for $P_C(t)$ and $P_C^{-1}(t)$ when $P_C(0) = 0$ are found 
by \eqref{inv-SG}
\begin{gather}
\hat{P}_1^C = \frac{1}{48}
\begin{bmatrix}
-1 & 1 & -1 \\
-1 & 1 & -1 \\
-1 & 1 & -1    
\end{bmatrix} , 
\hat{P}_2^C = \frac{1}{60}
\begin{bmatrix}
1 & -2 & 4 \\
2 & -4 & 8 \\
4 & -8 & 16    
\end{bmatrix}, \nonumber \\
\hat{P}_3^C = \frac{1}{240}
\begin{bmatrix}
-1 & 3 & -9 \\
-3 & 9 & -27 \\
-9 & 27 & -81    
\end{bmatrix} , 
\hat{P}_1^{-C} =  12
\begin{bmatrix}
-36 & 30 & -6 \\
-30 & 25 & -5 \\
-6  &  5 & -1    
\end{bmatrix} , \nonumber \\
\hat{P}_2^{-C} = 60
\begin{bmatrix}
9 & -12 & 3 \\
12 & -16 & 4 \\
3  &  -4 & 1    
\end{bmatrix} ,  
\hat{P}_3^{-C} = 60
\begin{bmatrix}
-4 & 6 & -2 \\
-6 & 9 & -3 \\
-2 & 3 & -1    
\end{bmatrix} . \nonumber 
\end{gather}
The orthogonality property $\hat{P}_j^{-C} \hat{P}_i^C = 0$ is satisfied for them when $i \ne j$.  
Therefore, at $t > 0$ we obtain
\begin{gather}
P_C^{-1}(t) P_C(t) = G(t) \sum_{j,i} \hat{P}_j^{-C} \hat{P}_i^C F_i(t) = G(t) \cdot \nonumber \\
\cdot \left( \hat{P}_1^{-C} \hat{P}_1^C F_1(t)+\hat{P}_2^{-C} \hat{P}_2^C F_2(t) + \hat{P}_3^{-C} \hat{P}_3^C F_3(t) \right) = I  \ \square .  \nonumber 
\end{gather}\\

\section{Extending Results to Multiple Eigenvalues in the Spectrum} 

In this section, we extend results obtained earlier to the case when the spectrum of matrix $A$ has multiple eigenvalues. 
This will allow us to obtain a closed-form estimate of the spectral properties of the system in the limiting case of resonant interaction, 
when the eigenvalues of the system converge.

We first extend the expressions \eqref{SD-gen-diff} and \eqref{subgram-general} for SDSE of Gramians of system \eqref{sys-0}.\\

\noindent \textbf{Theorem 6.} \textit{Let the eigenvalues $\lambda_i$ of matrix $A$ in \eqref{sys-0} have multiplicities $n_i$ $(n_1+\cdots+n_m =n)$.
Then SDSE of the infinite and finite Gramians $P(\infty)$ and $P(t)$ in \eqref{Lyap-algeb} and \eqref{Lyap-diff} are}
\begin{gather}
P(\infty) = \sum_{i=1}^m \sum_{k=1}^{n_i} \hat{A}_k^{(i)} \cdot BB^T \cdot (-I\lambda_i-A^T)^{-k} , \nonumber \\
P (t) = \sum_{i=1}^m \sum_{k=1}^{n_i} \hat{A}_k^{(i)}\cdot BB^T\cdot \Bigg[(- I \lambda_i-A^T)^{-k}  \nonumber \\
\left. - \sum_{l=1}^k (-I\lambda_i-A^T )^{-l} \frac{t^{k-l}}{(k-l)!} e^{(I\lambda_i+A^T )t} \right] , \label{mult-eig-0} \\
\hat{A}_k^{(i)} = 
\frac{1}{(n_i-k)!} \cdot \frac{d^{(n_i-k)}}{d\lambda^{(n_i-k)}} \left[(\lambda-\lambda_i)^{n_i} (I\lambda-A)^{-1} \right]_{\lambda=\lambda_i} , \nonumber
\end{gather}
\textit{where $\hat{A}_k^{(i)}$ are the matrix coefficients in the partial fractional decomposition of the resolvent.}\\
 
\noindent \textbf{Proof.} Up to notation, it is shown in \cite{DAN-2017} (see Corollary 5) that the spectral decomposition for $P(\infty)$ satisfies \eqref{Lyap-algeb}. The zero boundary condition $P_C (0) = 0$ can be checked directly from \eqref{mult-eig-0}. Note further that
\begin{gather}
A \hat{A}_{n_i}^{(i)} = \left[(-I\lambda+ A+I \lambda)(\lambda-\lambda_i)^{n_i} (I\lambda-A)^{-1} \right]_{\lambda=\lambda_i}  \nonumber \\
= - \left[(\lambda-\lambda_i)^{n_i} \right]_{\lambda=\lambda_i} + \lambda_i \hat{A}_{n_i}^{(i)} =  \lambda_i \hat{A}_{n_i}^{(i)} . \nonumber
\end{gather}
If $1\le k < n_i$ we have
\begin{gather}
 A \hat{A}_k^{(i)} = \frac{-1}{(n_i-k)!} \frac{d^{(n_i-k)}}{d\lambda^{(n_i-k)}} \left[(\lambda-\lambda_i)^{n_i} \right]_{\lambda=\lambda_i} \nonumber \\
+ \frac{1}{(n_i-k)!} \frac{d^{(n_i-k)}}{d\lambda^{(n_i-k)}} 
\left[ \lambda (\lambda-\lambda_i)^{n_i} (I\lambda-A)^{-1}  \right]_{\lambda=\lambda_i} = \nonumber \\
\lambda_i \hat{A}_k^{(i)} + \frac{n_i-k}{(n_i-k)!} \frac{d^{(n_i-k-1)}}{d\lambda^{(n_i-k-1)}}  
\left[(\lambda-\lambda_i)^{n_i} (I\lambda-A)^{-1} \right]_{\lambda=\lambda_i} \nonumber \\
= \lambda_i \hat{A}_k^{(i)} + \hat{A}_{k+1}^{(i)} . \nonumber
\end{gather}
Let us check the decomposition for $P(t)$. Note that the time-independent part of $P(t)$ satisfies 
the algebraic equation \eqref{Lyap-algeb}. Let us check that the exponential parts 
\begin{equation}
\tilde{P}_i (t) = - \sum_{k=1}^{n_i} \sum_{l=1}^k \hat{A}_k^{(i)} BB^T (-I\lambda_i-A^T )^{-l} \frac{t^{k-l}}{(k-l)!} e^{(I\lambda_i+A^T )t} \nonumber
\end{equation}
satisfy the differential Lyapunov equation \eqref{Lyap-diff} with a zero right-hand side. On the one hand,
\begin{gather}
-\frac{d\tilde{P}_i (t)}{dt} = - \tilde{P}_i (t) (I\lambda_i+A^T ) \nonumber \\
+ \sum_{k=1}^{n_i} \sum_{l=1}^{k-1} \hat{A}_k^{(i)} BB^T (-I\lambda_i-A^T )^{-l} 
\frac{t^{k-l-1}}{(k-l-1)!} \,e^{(I\lambda_i+A^T )t} . \nonumber
\end{gather}
On the other hand, using the expression for $A \hat{A}_k^{(i)}$ obtained earlier
\begin{gather}
A \tilde{P}_i (t) + \tilde{P}_i (t) A^T = \tilde{P}_i (t) (I\lambda_i+A^T )  \nonumber \\
- \sum_{k=1}^{n_i-1} \sum_{l=1}^k \hat{A}_{k+1}^{(i)} BB^T (-I\lambda_i-A^T )^{-l} 
\frac{t^{k-l}}{(k-l)!} \,e^{(I\lambda_i+A^T )t} . \nonumber
\end{gather}
After replacing here the index $\tilde{k}=k+1$, we see that the sum of both expressions above gives zero. Therefore, the time-dependent part of $P (t)$, i.e., $\sum_i \tilde{P}_i (t)$, 
satisfies the Lyapunov differential equation \eqref{Lyap-diff} with a zero right-hand side. $\blacksquare$\\

To obtain the SDSE of $P_C$ and $P^{-1}_C$ for the system in the controllability canonical form, we will use the representation 
of matrix $A$ through generalized eigenvectors. Let $m$ eigenvalues in spectrum $\sigma(A) = \{\lambda_1,\lambda_2,\dots,\lambda_m\}$ have multiplicities $n_1, n_2,\dots, n_m$, 
respectively ($n= \sum_i n_i$). Then matrix $A$ can be represented in a Jordan normal form
\begin{equation}
A = M J M^{-1}, \ J = J_1 \oplus J_2 \oplus \cdots \oplus J_m, \ \nonumber 
\end{equation}
where $J_i$ is a Jordan block corresponding to $\lambda_i$, $\oplus$ is a direct sum of matrices, 
and matrices $M$ and $M^{-1}$ are composed of generalized right and left eigenvectors, respectively 
\begin{gather}
M = [x^{(1)}_1, \dots, x^{(1)}_{n_1}; x^{(2)}_1, \dots, x^{(2)}_{n_2}; \dots; x^{(m)}_1, \dots, x^{(m)}_{n_m}], \nonumber \\
M^{-1} = [y^{(1)}_1, \dots, y^{(1)}_{n_1}; y^{(2)}_1, \dots, y^{(2)}_{n_2}; \dots; y^{(m)}_1, \dots, y^{(m)}_{n_m}]^T. \nonumber
\end{gather}
We assume them to be normalized as $x^{(i)}_k (y^{(j)}_l)^T = \delta_{ij} \delta_{kl}$.
Each Jordan chain of eigenvectors 
we denote as
\begin{gather}
M_i = [x^{(i)}_1, \dots, x^{(i)}_{n_i}] \ \text{and} \ M^{(-1)}_i = [y^{(i)}_1, \dots, y^{(i)}_{n_i}]^T \ \ \text{so that} \nonumber \\
M^{(-1)}_i M_j =\delta_{ij} I_{n_i\times n_i}, \ \ M_i M^{(-1)}_i = \sum^{n_i}_{k=1} x^{(i)}_k (y^{(i)}_k)^T \, . \nonumber
\end{gather} 
The eigenvectors in the $i$-th Jordan chain can be obtained successively from
\begin{gather}
(A-\lambda_i I) x^{(i)}_1 = 0, \ (A-\lambda_i I)x^{(i)}_2 = x^{(i)}_1, \dots, \nonumber \\
 \ (A-\lambda_i I)x^{(i)}_{n_i} = x^{(i)}_{n_i-1}. \label{ev-J-chain}
\end{gather}
If the system is given in the controllability canonical form, 
then from \eqref{ev-J-chain} one can calculate all eigenvectors explicitly as functions of $\lambda_i$:
\begin{gather}
x^{(i)}_1 =\frac{1}{\lambda_i} [ \lambda_i, \lambda_i^2, \cdots, \lambda_i^n]^T, \nonumber \\
x^{(i)}_2 =\frac{1}{\lambda^2_i} [ \lambda_i, 2 \lambda_i^2, \cdots, k \lambda^k_i, \cdots, n \lambda_i^{n}]^T , \nonumber \\
x^{(i)}_3 =\frac{1}{\lambda^3_i} [ \lambda_i, 2 \lambda_i^2,  \cdots, \left(1+\frac{n(n-1)}{2}\right) \lambda_i^{n}]^T , \ \text{etc.} \label{J-chain-c} 
\end{gather}
Now we are ready to extend SDSE in \eqref{inv-SG} in Theorem 5 to the case of multiple eigenvalues in the spectrum of $A_C$.\\
  
\noindent \textbf{Theorem 7.} \textit{The SDSE of the infinite Gramian $P_C(\infty)$ and its inverse $P^{-1}_C(\infty)$
of system \eqref{sys-c}, satisfying Lyapunov and Riccati algebraic equations
\eqref{Lyap-algeb-c} and \eqref{Ric-algeb-c}, are given by}
\begin{gather}
P_{C}(\infty) = \sum^m_{i=1} \hat{P}_i^C, \ \ P^{-1}_{C}(\infty) = \sum^m_{j=1} \hat{P}_j^{-C} , \nonumber \\ 
\hat{P}_i^C = M_i \mathcal{H}_i \mathcal{T}_i^{-1} M^T_i  \mathcal {J}, \nonumber \\
\hat{P}_j^{-C} = \mathcal{J} (M^{(-1)}_j)^T \mathcal{T}_j \mathcal{H}^{-1}_j M^{(-1)}_j, \label{SDSE-mult-c} \\
\text{where} \ \ \mathcal{T}_i =
\begin{bmatrix}
c^T x^{(i)}_1 & 0 & \cdots & 0 \\
c^T x^{(i)}_2 & c^T x^{(i)}_1  & \ddots &\vdots \\
\vdots & \ddots & \ddots & 0 \\
c^T x^{(i)}_{n_i} & \cdots & c^T x^{(i)}_2 & c^T x^{(i)}_1
\end{bmatrix},   \label{matrix-T} \\
c^T = a^T((-1)^nI+ \mathcal{J}), \ \  a^T =[a_0, a_1, \dots, a_{n-1}] , \nonumber \
\end{gather}
\begin{equation}
\mathcal{H}_i =
\begin{bmatrix}
e^T_n y^{(i)}_1 & e^T_n y^{(i)}_2 & \cdots & e^T_n y^{(i)}_{n_i} \\
e^T_n y^{(i)}_2 & & \iddots  & 0 \\
\vdots & \iddots & \iddots &  \vdots \\
e^T_n y^{(i)}_{n_i} & 0 & \cdots &  0
\end{bmatrix} , \label{matrix-H} 
\end{equation}
\textit{$\mathcal{H}_i$ and $\mathcal{T}_i$ are the upper Hankel and lower Toeplitz matrices, respectively, and eigenparts in \eqref{SDSE-mult-c} satisfy the normalization condition $\hat{P}_i^C \hat{P}_j^{-C} =\delta_{ij} M_i M^{(-1)}_i$.}\\

\noindent \textbf{Proof.} 
 According to Theorem 6, the Lyapunov algebraic equation for system \eqref{sys-c} is satisfied if
\begin{equation}
\hat{P}_i^C =\left[ \hat{A}_1^{(i)} e_n,  \cdots, \hat{A}_{n_i}^{(i)} e_n \right]
\begin{bmatrix}
e^T_n (-\lambda_i I - A^T_C)^{-1} \\
e^T_n (-\lambda_i I - A^T_C)^{-2} \\
\cdots \\
e^T_n (-\lambda_i I - A^T_C)^{-n_i} \\
\end{bmatrix}. \label{prod-vect}
\end{equation}
Consider the first vector in \eqref{prod-vect}. 
Let us prove that it is equal to $M_i \mathcal{H}_i$. To do this, we prove that
\begin{equation}
\hat{A}_k^{(i)} = \sum^{n_i+1-k}_{l=1} x^{(i)}_l (y^{(i)}_{k-1+l})^T . \label{Ak-matrix}
\end{equation}
Matrices $\hat{A}_k^{(i)}$ are defined as coefficients in the partial fractional decomposition 
\begin{equation}
(I\lambda-A)^{-1} = \sum^m_{i=1} \sum^{n_i}_{k=1} \frac{\hat{A}_k^{(i)}}{(\lambda-\lambda_i)^k}. \nonumber
\end{equation}
Let us prove that \eqref{Ak-matrix} satisfies this definition, i.e.,
\begin{equation}
(I\lambda-A)^{-1} = \sum^m_{i=1} \sum^{n_i}_{k=1} \sum^{n_i+1-k}_{l=1} \frac{x^{(i)}_l (y^{(i)}_{k-1+l})^T}{(\lambda-\lambda_i)^k} . \label{eq-1}
\end{equation}
From the Jordan normal form
\begin{equation}
A = M J M^{-1} = \sum^m_{i=1} \left( \sum^{n_i}_{k=1} \lambda_i x^{(i)}_k (y^{(i)}_k)^T +\sum^{n_i-1}_{k=1} x^{(i)}_k (y^{(i)}_{k+1})^T \right) \nonumber
\end{equation}
it follows that
\begin{gather}
I \lambda-A = \lambda M M^{-1} - M J M^{-1} \nonumber \\
= \sum^m_{i=1} \left( \sum^{n_i}_{k=1} (\lambda-\lambda_i) x^{(i)}_k (y^{(i)}_k)^T - \sum^{n_i-1}_{k=1} x^{(i)}_k (y^{(i)}_{k+1})^T \right) . \label{eq-2}
\end{gather}
Let us multiply \eqref{eq-1} and \eqref{eq-2} using orthogonality of eigenvectors 
\begin{gather}
 \sum^m_{i,i'=1} \sum^{n_i}_{k=1} \sum^{n_i+1-k}_{l=1} \frac{x^{(i)}_l (y^{(i)}_{k-1+l})^T}{(\lambda-\lambda_i)^k} \times \nonumber \\
 \times \left[  \sum^{n_i}_{k'=1} (\lambda-\lambda_{i'}) x^{(i')}_{k'} (y^{(i')}_{k'})^T - \sum^{n_i-1}_{k'=1} x^{(i')}_{k'} (y^{(i')}_{k'+1})^T \right] \nonumber \\
 =  \sum_{i} \left[ \sum^{n_i}_{k=1} \sum^{n_i+1-k}_{l=1} \frac{x^{(i)}_l (y^{(i)}_{k-1+l})^T}{(\lambda-\lambda_i)^{k-1}}
 - \sum^{n_i-1}_{k=1} \sum^{n_i-k}_{l=1} \frac{x^{(i)}_l (y^{(i)}_{k+l})^T}{(\lambda-\lambda_i)^{k}} \right] \nonumber \\
= \sum_i \left[  
 \sum^{n_i}_{l=1} x^{(i)}_l (y^{(i)}_l)^T 
 + \sum^{n_i - 1}_{\tilde{k}=1} \sum^{n_i-\tilde{k}}_{l=1} \frac{x^{(i)}_l (y^{(i)}_{\tilde{k}+l})^T}{(\lambda-\lambda_i)^{\tilde{k}}} \right. \nonumber \\
 \left. - \sum^{n_i-1}_{k=1} \sum^{n_i - k}_{l=1} \frac{x^{(i)}_l (y^{(i)}_{k+l})^T}{(\lambda-\lambda_i)^k} \right] 
 =  \sum_i \sum^{n_i}_{l=1} x^{(i)}_l (y^{(i)}_l)^T = I . \nonumber 
\end{gather}
The obtained identity proves \eqref{Ak-matrix}; therefore, we obtain
\begin{equation}
\left[ \hat{A}_1^{(i)} e_n,  \cdots, \hat{A}_{n_i}^{(i)} e_n \right] = M_i \mathcal{H}_i . \label{1-vector}
\end{equation}

Consider the second vector in \eqref{prod-vect}. Let $a =[a_0, a_1, \dots, a_{n-1}]$ be the vector of coefficients of the characteristic polynomial. 
Then $ A_C = A_0 - e_n a^T$, where $A_0$ is a matrix of zeros except units over the diagonal, 
and $e_n$ is the last column of the unit matrix. According to \eqref{ev-J-chain}, for the first eigenvector corresponding to $\lambda_i$ we have
\begin{gather}
(A_0 - e_n a^T -\lambda_i I ) x^{(i)}_1 = 0  \ \ \Rightarrow \nonumber \\
(\mathcal{J}A_0\mathcal{J} - \mathcal{J}e_n a^T\mathcal{J} -\lambda_i I ) \mathcal{J} x^{(i)}_1 = 0 \ \ \Rightarrow \nonumber \\
(- A_0 + e_n a^T - e_n a^T -(-1)^n e_n a^T\mathcal{J} -\lambda_i I ) \mathcal{J} x^{(i)}_1 = 0 \ \ \Rightarrow \nonumber \\
( - A_C -\lambda_i I)  \mathcal{J} x^{(i)}_1 = e_n c^T x^{(i)}_1  \ \ \Rightarrow \nonumber \\
(x^{(i)}_1)^T \mathcal{J} = c^T x^{(i)}_1 e^T_n (-\lambda_i I - A^T_C)^{-1} . \nonumber
\end{gather}
Similarly, 
\begin{gather}
(A_C -\lambda_i I ) x^{(i)}_2 = x^{(i)}_1  \ \ \Rightarrow \nonumber \\
( - A_C -\lambda_i I)  \mathcal{J} x^{(i)}_2 = e_n c^T x^{(i)}_2 + \mathcal{J} x^{(i)}_1 \ \ \Rightarrow \nonumber \\
(x^{(i)}_2)^T \mathcal{J} = c^T x^{(i)}_2 e^T_n (-\lambda_i I - A^T_C)^{-1} \nonumber \\
+ c^T x^{(i)}_1 e^T_n (-\lambda_i I - A^T_C)^{-2} . \nonumber
\end{gather}
Continuing the same operation up to $x^{i}_{n_i}$, we obtain 
\begin{equation}
M_i^T \mathcal{J} = \mathcal{T}_i  
\begin{bmatrix}
e^T_n (-\lambda_i I - A^T_C)^{-1} \\
e^T_n (-\lambda_i I - A^T_C)^{-2} \\
\cdots \\
e^T_n (-\lambda_i I - A^T_C)^{-n_i} \\
\end{bmatrix}, \label{2-vector}
\end{equation}

From \eqref{1-vector} and \eqref{2-vector} we obtain an expression for $\hat{P}_i^C$.
The expression for $\hat{P}_j^{-C}$ is chosen to satisfy the normalization condition. 
In this case
\begin{gather}
P_C P^{-1}_C =\sum^m_{i,j=1} \hat{P}_i^C \hat{P}_j^{-C} = \sum^m_{i=1} M_i M^{(-1)}_i \nonumber \\
= \sum^m_{i=1} \sum^{n_i}_{k=1} x^{(i)}_k (y^{(i)}_k)^T = M M^{-1} = I \, .  \nonumber
\end{gather}
From here it follows that algebraic Riccati equation \eqref{Ric-algeb-c} is also satisfied for $P^{-1}_C$ $\blacksquare$.\\

\noindent \textbf{Remark 10.} Although the numerical computation of the Jordan canonical form is an ill-posed problem and the eigenvectors may be chosen in more than one way \cite{Antoulas-2005}, expressions \eqref{SDSE-mult-c} are uniquely determined, since the derivation uses the unambiguous expressions of Theorem 6 and the recursive relation \eqref{ev-J-chain}, which is valid for any choice of Jordan chain of eigenvectors. In fact, any convenient algorithm can be used to compute the generalized eigenvectors, as is done in \eqref{J-chain-c}.\\

\noindent \textbf{Example 5} illustrates and verifies the formulas of Theorem 7.
Consider a system with eigenvalue $\lambda_1=1$ of multiplicity 2 and $\lambda_2=2$ of multiplicity 3 in the controllability canonical form. We have 
\begin{gather}
N(s) = s^5 - 8s^4 + 25 s^3 - 38s^2 + 28s - 8, \nonumber \\
a^T =[-8,28,-38,25,-8], \ c^T=[16,0,79,0,16] . \nonumber
\end{gather}
According to \eqref{J-chain-c} we obtain generalized eigenvectors
\begin{equation}
M=[M_1,M_2] =
\begin{bmatrix}
1 & 1 & 1  & 0.5 & 0.25 \\
1 & 2 & 2  & 2 & 1 \\
1 & 3 & 4  & 6 & 4 \\
1 & 4 & 8  & 16 &14 \\
1 & 5 & 16 & 40 & 44
\end{bmatrix}, \nonumber 
\end{equation}
\begin{equation}
M^{-1}=
\begin{bmatrix}
M^{-1}_1 \\
M^{-1}_2 
\end{bmatrix} =
\begin{bmatrix}
-8 & 28 & -30  & 13 & -2 \\
-8 & 20 & -18  & 7 & -1 \\
22 & -62.5 & 63 & -26.5 & 4 \\
-12 & 35 & -36.5  & 16 & -2.5 \\
 4 & -12 & 13 & -6 & 1
\end{bmatrix} . \nonumber
\end{equation}
According to \eqref{matrix-T} and \eqref{matrix-H} we obtain
\begin{gather}
\mathcal{T}_1 =
\begin{bmatrix}
108 & 0 \\
324 & 108 
\end{bmatrix}, \ \mathcal{H}_1
\begin{bmatrix}
-2 & -1 \\
-1 & 0 
\end{bmatrix}, \nonumber \\
\mathcal{T}_2 =
\begin{bmatrix}
576 & 0 &0 \\
1104 & 576 & 0 \\
1012 & 1104 & 576  
\end{bmatrix}, \ \mathcal{H}_2
\begin{bmatrix}
4 & -2.5 & 1 \\
-2.5 & 1& 0 \\
1 & 0 & 0
\end{bmatrix}. \nonumber 
\end{gather}
According to \eqref{SDSE-mult-c} we obtain the symmetrized sub-Gramians
\begin{gather}
\tilde{P}^C_1 =
\left\{ \hat{P}^C_1 \right\}_H = \frac{1}{108}
\begin{bmatrix}
1 &  0 & 3 &  0 & 5 \\
0 & -3 & 0 & -5 & 0 \\
3 &  0 & 5 &  0 & 7 \\
0 & -5 & 0 & -7 & 0 \\
5 &  0 & 7 &  0 & 9
\end{bmatrix}, \  \nonumber \\
\tilde{P}^C_2 = \frac{1}{128\cdot108}
\begin{bmatrix}
-169 &  0 & -372 &  0 & -656 \\
0      & 372 & 0   & 656 & 0 \\
-372 &  0   & -656 &  0 & -832 \\
0      & 656 & 0    & 832 & 0 \\
-656 &  0    & -832 &  0 & -2304 
\end{bmatrix}, \  \nonumber
\end{gather}
Their sum 
$\tilde{P}^C_1+\tilde{P}^C_2$
provides an exact solution of Lyapunov equation \eqref{Lyap-algeb-c}
\begin{equation}
P_C = 
\begin{bmatrix}
-41 &  0  & 12 &  0 & -16 \\
  0 & -12 & 0 & 16 & 0 \\
12 &  0 & -16 &  0 & 64 \\
0   & 16 & 0   & -64 & 0 \\
-16 &  0 & 64 &  0  & 1152
\end{bmatrix} . \  \nonumber 
\end{equation}
According to \eqref{SDSE-mult-c} we obtain the symmetrized eigenparts for the inverse of the Gramian
\begin{gather}
\tilde{P}^{-C}_1 = 108 
\begin{bmatrix}
192 &  0 & 528 &  0 & 32 \\
0 & -1520 & 0 & -596 & 0 \\
528 &  0 & 1404 &  0 & 84 \\
0 & -596 & 0 & -231 & 0 \\
32 &  0 & 84 &  0 & 5
\end{bmatrix}, \  \nonumber \\
\tilde{P}^{-C}_2 = 4
\begin{bmatrix}
-5296 &  0 & -14356 &  0 & -868 \\
0      & 40608 & 0   & 15984 & 0 \\
-14356 &  0   & -38275 &  0 & -2287 \\
0      & 15984 & 0    & 6156 & 0 \\
-868 &  0    & -2287 &  0 & -139 
\end{bmatrix}. \nonumber
\end{gather}
Their sum 
$\tilde{P}^{-C}_1+\tilde{P}^{-C}_2$
provides an inverse of $P_C$ 
and the solution of the Riccati equation \eqref{Ric-algeb-c} in closed-form $\square$.\\

\section{Conclusion}

In this paper, we obtain SDSE of the controllability Gramian and its inverse for a continuous LTI dynamical system in the canonical controllability form. These decompositions represent closed-form solutions of the algebraic and differential Lyapunov and Riccati equations. The spectral components in these decompositions allow for a quantitative estimation of the controllability, observability, and stability properties of the system in connection with its individual eigenmodes and their pairwise combinations, and relate these properties to the corresponding devices of the system based on modal analysis. 

From a mathematical point of view, the symmetrized eigenparts in Gramian expansions are uniquely determined as amplitudes at different eigenmodes in the solution of the Lyapunov differential equation (see Theorem 2). Conversely, the eigenparts in the expansion of the inverse of Gramian can be uniquely determined based on the orthogonality condition between the eigenparts of the Gramian and its inverse (see Theorem 4). 
For the inverse of the finite Gramian, there is a certain non-uniqueness associated with the optimal choice of the normalization matrix $G(t)$ in Theorem 5, which requires further investigation. From an engineering point of view, as proven in Properties 1 and 4, the obtained symmetrized spectral components have a physical interpretation as measurable quantities in the minimum energy control problem. Therefore, they are unambiguously defined. 

Three main results are obtained in the paper: First, the closed-form SDSE of the inverse of controllability Gramian in Theorems 4 and 5 is obtained for the first time to the best of our knowledge. 
This decomposition can quantitatively characterize the influence of individual eigenmodes in the system and the specific devices associated with them on the minimum control energy.
Second, the SDSE of Gramians obtained earlier in the literature are extended to solutions of differential Lyapunov and Riccati equations with arbitrary initial conditions (Theorems 2 and 5, Remarks 3 and 4), which allows for the estimation of system spectral properties over an arbitrary time interval and their prediction at future moments. 
Finally, the SDSE of the Gramian and its inverse are generalized to the case of multiple eigenvalues in the spectrum of the dynamics matrix (Theorems 6 and 7), which allows us to estimate in closed form the effects of resonant interactions with the eigenmodes of the system.

We expect that the additional information obtained from SDSEs of Gramians and their inverses can improve the accuracy of algorithms in solving various practical problems such as minimum energy control \cite{Klamka-2018}, design in aerospace engineering \cite{Gupta-2020}, tuning power system stabilizers \cite{Ghosh-2020}, model order reduction \cite{Benner-2021}, \cite{Casadei-2020}, design of structural metrics for complex networks \cite{Summers-2016}, stability analysis of inter-area oscillations in power systems \cite{Isk-2023}, optimal actuator and sensor placement \cite{Dilip-2019}, and traffic network optimization \cite{Bianchin-2020}.

\section*{Acknowledgment}
This work was supported by the Russian Science Foundation, project no. 25-29-20158.
We thank M.V. Morozov for valuable remarks.


\end{document}